\documentclass[11pt,review]{elsarticle}
\journal{}
\usepackage{lineno,hyperref}
\usepackage[margin=25mm]{geometry}
\usepackage{amsmath}
\usepackage{amsfonts}
\usepackage{amssymb}

\usepackage{verbatim}
\usepackage{multirow}
\usepackage{natbib}
\usepackage{pbox}
\immediate\write18{texcount -tex -sum  \jobname.tex > \jobname.wordcount.tex}
\usepackage[parfill]{parskip}
\usepackage{hyperref}
\usepackage{algorithm}
\usepackage[noend]{algpseudocode}
\usepackage{pgfplots}
\usepackage{tikz}
\usetikzlibrary{calc,positioning,shapes,fit}
\pgfplotsset{compat=1.17} 

\usepackage{setspace}
\tikzset{pblock/.style = {rectangle split,
                      rectangle split parts=2, very thick,draw=black!50, top
                      color=white,bottom color=black!20, align=center}}
 
 \tikzset{ldgrayblock/.style = {rectangle split,
                      rectangle split parts=2, very thick,draw,rectangle split part fill={white,black!20}, align=center}}

\tikzset{twoblock/.style = {rectangle split,
                      rectangle split parts=2, very thick,draw, align=center}}        
                      
\tikzset{longtwoblock/.style = {rectangle split,
                      rectangle split parts=2, very thick, align=center, draw, text width=12cm}}                       
                      
\tikzstyle{process} = [rectangle, minimum width=4cm, minimum height=1cm, inner sep=3mm, text centered,  very thick, align=center,draw=black]

\tikzstyle{crect} = [rectangle, rounded corners, minimum width=2cm, minimum height=1cm, inner sep=3mm, text centered,  very thick, align=center,draw=black]
 
\tikzstyle{decision} = [diamond, aspect=1, text centered,  very thick, align=center,  minimum width=1mm, minimum height=1mm, draw=black]

\tikzstyle{circmarker}=[circle,  text centered,  very thick, align=center,  minimum width=2mm,draw=black,fill=white,inner sep=0.5mm]

\tikzstyle{circstate}=[circle,  text centered,  very thick, align=center,  minimum width=1mm,draw=black,fill=white,inner sep=1mm]

\tikzstyle{circH}=[circle,  text centered,  very thick, align=center,  minimum width=3mm,draw=black,fill=white,inner sep=1mm]

 \tikzset{every fit/.style={text width=20cm, inner sep=3mm}}
\hypersetup{ %
	pdftitle={L2P},
	pdfkeywords={},
	pdfborder=0 0 0,
	pdfpagemode=UseNone,
	colorlinks=true,
	linkcolor=blue, 
	citecolor=blue, 
	filecolor=blue, 
	urlcolor=blue, 
	pdfview=FitH,
	pdfauthor={Anonymous},
}

\DeclareMathOperator*{\argmax}{arg\,max}
\DeclareMathOperator{\CLIP}{CLIP}
\DeclareMathOperator{\Iter}{Iter}
\DeclareMathOperator{\iter}{iter}
\DeclareMathOperator{\newBest}{newBest}
\DeclareMathOperator{\true}{true}
\DeclareMathOperator{\false}{false}
\DeclareMathOperator{\Uniform}{Uniform}
\DeclareMathOperator{\prob}{prob}
\def\RL{\mathrm{RL}}

\def\Temp{\mathrm{Temp}}
\def\best{\mathrm{best}}
\def\cooling{\mathrm{cooling}}
\def\MAX{\mathrm{MAX}}
\def\TERMINAL{\mathrm{TERMINAL}}
\def\start{\mathrm{start}}

\def\old{\mathrm{old}}
\def\gap{\mathrm{gap}}

\usepackage{array}
\usepackage{makecell}

\begin{document}

\begin{frontmatter}
\title{Learning to Schedule Heuristics for the Simultaneous Stochastic Optimization of Mining Complexes}
\author{Yassine Yaakoubi$^{1,2}$, Roussos Dimitrakopoulos$^{1,2}$  \\
        \small $^{1}$COSMO --- Stochastic Mine Planning Laboratory, McGill University \\
        \small $^{2}$GERAD --- Research Group in Decision Analysis\\
}
\begin{abstract}
The simultaneous stochastic optimization of mining complexes (SSOMC) is a large-scale stochastic combinatorial optimization problem that simultaneously manages the extraction of materials from multiple mines and their processing using interconnected facilities to generate a set of final products, while taking into account material supply (geological) uncertainty to manage the associated risk. Although simulated annealing has been shown to outperform comparing methods for solving the SSOMC, early performance might dominate recent performance in that a combination of the heuristics' performance is used to determine which perturbations to apply. This work proposes a data-driven framework for heuristic scheduling in a fully self-managed hyper-heuristic to solve the SSOMC. The proposed learn-to-perturb (L2P) hyper-heuristic is a multi-neighborhood simulated annealing algorithm. The L2P selects the heuristic (perturbation) to be applied in a self-adaptive manner using reinforcement learning to efficiently explore which local search is best suited for a particular search point. Several state-of-the-art agents have been incorporated into L2P to better adapt the search and guide it towards better solutions. By learning from data describing the performance of the heuristics, a problem-specific ordering of heuristics that collectively finds better solutions faster is obtained. L2P is tested on several real-world mining complexes, with an emphasis on efficiency, robustness, and generalization capacity. Results show a reduction in the number of iterations by $30$--$50\%$ and in the computational time by $30$--$45\%$.
\end{abstract}

\end{frontmatter}

\section{Introduction}

Solving a large-scale industrial mining complex problem is critical for decision-makers in the mining industry. An industrial mining complex is an engineering system that manages the extraction of materials from multiple mines, their processing in interconnected facilities and multiple waste sites to generate a set of products delivered to customers. These related operations require a strategic mine plan that maximizes the net present value (NPV) of mining operations while meeting various complex constraints. Mine scheduling optimization is critical as sub-optimal scheduling can adversely affect profit.
Conventional mine planning methods break down the optimization into distinct steps to be solved sequentially. They typically optimize mine production scheduling for individual mines (defining the sequence of block extraction), then the downstream (destination) policies. This decomposition leads to suboptimal results because it does not account for the various interactions between processing, transportation, and destination policies contributing to the value of mined products. Conventional methods tend to ignore material supply uncertainty to reduce the model size to allow the use of commercial solvers, which is known to be the main cause of failure to meet production targets and forecasts~\cite{godoy2003,dimitrakopoulos2011}. Conventional methods also tend to aggregate blocks~\cite{ramazan2007Tree, froyland2009, stone2018}, leading to homogenization of spatial distribution of metal within an aggregate.

A new framework termed simultaneous stochastic optimization of mining complexes (SSOMC) has been proposed~\cite{montiel2015,montiel2017,goodfellow2016,goodfellow2017} to optimize the system by incorporating all components of the mining complex into a single mathematical model. The SSOMC model maximizes the NPV of the mining complex and includes geological uncertainty to manage the associated risk. This type of uncertainty is modeled using a set of geostatistically (stochastically) simulated scenarios of all relevant mineral deposit attributes~\cite{goovaerts1997, remy2009, david2012, rossi2013}. Unlike conventional methods that consider the economic value of blocks, the SSOMC considers the simultaneous extraction of material from multiple mines and optimizes the economic value of the final products and all intermediate interactions. Since exact methods are impractical to use for SSOMC with more than a few thousand mining blocks, researchers have employed metaheuristics and hyper-heuristics to address the SSOMC (as well as mine production scheduling). Most notably, a simulated annealing approach proposed by~\citet{goodfellow2016} outperforms all comparing methods. However, early performance may dominate recent performance in that the score function only uses a combination of past performance to sample the perturbations to be applied and solution time remains an issue.

Machine learning (ML) has shown promise for addressing small and large-scale operations research (OR) tasks~\citep{Yaakoubi2019, gasse2019, Yaakoubi2020, gupta2020, Yaakoubi2021}. Recent research even uses ML to guide OR algorithms to solve combinatorial optimization problems~\citep{zarpellon2020, chmiela2021, lipets2021}. Given the complexity of these problems, state-of-the-art algorithms rely on specific heuristics to make decisions otherwise computationally expensive or mathematically intractable. Therefore, ML emerges as a suitable predictive tool for making such heuristic scheduling decisions in a more reasoned and optimized manner~\citep{bengio2020}. In this manuscript, this line of work is followed by proposing, to the best of our knowledge, the first data-driven framework for heuristic scheduling in a fully self-managed hyper-heuristic-based solver. Indeed, in most of the data-driven AI for OR literature~\cite{Khalil2017,bengio2020}, the authors assume structural knowledge of the scheduling problem (e.g., optimal solutions or strong branching decisions). This assumption allows to apply ML to provide fast approximations of these decisions. Limited research has been conducted to solve large-scale stochastic combinatorial optimization problems such as the SSOMC where no structural knowledge pre-exists. In this case, ML models can be deployed to discover the structure of the problem. This is called \textit{policy learning} and is an application of reinforcement learning (RL).

The proposed self-learning hyper-heuristic, L2P (learn-to-perturb), improves the search strategy by using past experience to better guide the exploration, thus accelerating the process. This hyper-heuristic uses a multi-neighborhood simulated annealing algorithm in conjunction with RL. By defining a neighborhood structure, the RL agent learns to guide the search through the solution space landscape to find better solutions during the optimization process.
This paper aims to contribute the following: (1) propose a self-managed hyper-heuristic for solving the SSOMC based on a multi-neighborhood simulated annealing algorithm with adaptive neighborhood search; (2) use RL to adapt the search and guide it towards better solutions; (3) test on several mining complexes with emphasis on efficiency, robustness, and generalization capacity. To test the proposed methodology, three policy-based agents are proposed: advantage actor-critic (A2C)~\citep{a2c}, proximal policy optimization (PPO)~\citep{ppo}, and soft actor-critic (SAC)~\citep{sac}. Results demonstrate the effectiveness of L2P on a series of real-world mining complexes, where the proposed solution process takes only minutes, whereas CPLEX solves the linear relaxation in days or weeks. Using RL further reduces the computational time by $30$--$45\%$.

The rest of the paper is structured as follows: Section~\ref{sec:relWork} reviews the literature on metaheuristics and hyper-heuristics in the context of the SSOMC. Section~\ref{sec:math_form} presents an overview of the mathematical formulation of the SSOMC. Section~\ref{sec:method} details the L2P solution approach. Then, Section~\ref{sec:comp_exp} presents computational experiments on small- and large-scale instances of the SSOMC, along with comparisons for L2P variants (either using a standard score function or using various RL agents to guide the search), followed by conclusions in Section~\ref{sec:conclusion}.

\section{Literature review}
\label{sec:relWork}

\subsection{Mathematical programming for the simultaneous stochastic optimization of mining complexes}

Strategic or long-term mine planning aims to optimize when to extract material from mineral deposits and where to send it under various technical constraints in order to maximize the economic value of extraction. Such planning uses simulated mineral deposits that are discretized into a three-dimensional block model where each block is assigned values of spatially distributed attributes of interest (e.g., grade). Since this problem is of utmost importance to the mining industry, many papers address it~\cite{Hoerger1999, whittle2007, whittle2010, zuckerberg2011, topal2012, stone2018}. Due to the large size of orebody models and the complexity of constraints within, the industry practice for many years has been to divide the strategic mine planning into several disconnected subproblems and solve them sequentially using approximate methods with strong assumptions, ultimately leading to sub-optimal solutions. Previous research highlights not only the importance of linking the multiple components of the mineral value chain using a single simultaneous optimization framework but also the incorporation of geological uncertainty in a stochastic formulation, which increases the net present value (NPV) by 15-30\%~\cite{lamghari2012, ramazan2013, leite2014, goodfellow2016}, compared to assuming a deterministic input of grade distribution.

To overcome the above limitations, the SSOMC~\cite{montiel2015,montiel2017,goodfellow2016,goodfellow2017} accounts for geological uncertainty using two-stage stochastic programming models and links all components of a mining complex and can model non-linear interactions at each destination. Unlike previous work that considers the economic value of a block (due to dividing the problem into several sub-problems to be solved sequentially), the SSOMC considers revenues from final products and all intermediate interactions (blending, various costs along the value chain, any non-linear interactions such as in processors and stockpiles). As such, the SSOMC optimizes not only the block extraction sequence (production scheduling) but also the decision policies, downstream optimization (and potentially transportation), in a single model, capitalizing on the synergies and interactions between the different components.

Since the SSOMC integrates all components of the value chain into one model (from mines to products), the problem is too complex to be solved by exact methods. Thus, it becomes important to use heuristic methods, and mainly metaheuristics have been used to solve strategic mine planning in general, and the SSOMC in particular. Metaheuristics combine (problem-specific) heuristics into a more general framework. Thus, although they do not guarantee a global optimum in the presence of imperfect information and resource constraints, they may find a near-optimal solution. Fundamental metaheuristics include particle swarm optimization~\citep{PSO} inspired by birds flocking, ant colony optimization~\citep{ACO} inspired by the interactive behavior of ants and simulated annealing~\citep{simulated2} inspired by annealing process in metallurgy. Since these metaheuristics have been successfully applied to open-pit mine scheduling (i.e., production scheduling), such as Tabu search~\cite{lamghari2012} and variable neighborhood descent algorithm~\citep{lamghari2014}, they were then extended to the SSOMC.

Two main approaches have been developed to address the SSOMC using metaheuristics. The first approach~\cite{montiel2015,montiel2017} proposes to simultaneously optimize the mining, processing, and transportation schedules by maximizing the NPV of the mining complex and minimizing the deviations from the ore and waste production targets. The proposed two-stage stochastic integer program (SIP) formulation is solved using an efficient meta-heuristic, which considers three sequentially executed perturbation strategies: block-based perturbations (modifying the extraction period and destination of each block), operating alternatives-based perturbations (allowing for better control of geometallurgical variables), and transpiration system perturbations (modifying the proportions of materials through the available logistics system to deliver materials from one destination to another).

In the second approach~\cite{goodfellow2016,goodfellow2017}, the authors propose a formulation of a two-stage nonlinear SIP optimizing three decision variables. In addition to the scenario-independent block extraction sequence and scenario-dependent processing stream decisions (proportions sent from one location to another), the authors use k-means~\cite{macqueen1967,lloyd1982} to group blocks with similar characteristics (attributes) and send them to the same destination. This clustering allows destination policies to be dynamically changed during optimization and reduces the number of processing decisions that need to be made, in contrast to the blockwise destination SIP optimization in~\citet{montiel2015} and~\citet{montiel2017}. \citet{goodfellow2016,goodfellow2017} compare three metaheuristic algorithms, including simulated annealing, particle swarm optimization, and differential evolution, achieving an NPV up to 13\% higher than standard deterministic industrial software. The framework is particularly flexible, allowing non-linear calculation of attributes (blending, grade recovery function, geometallurgical responses) at any stage of the mining complex.

A relatively recent development in metaheuristic optimization is that of hyper-heuristics~\citep{burke2013hyper}. Hyper-heuristics have become particularly interesting since many existing approaches rely on metaheuristics, and no single approach works for various types of mining complexes. Hyper-heuristics typically combine low-level heuristics to perturb the current solution, so they operate at a higher level of abstraction than a metaheuristic~\citep{burke2013hyper}. Hyper-heuristics have been successfully employed for mine scheduling (but not yet for SSOMC). For example,~\citet{lamghari2018} use 27 low-level heuristics to compare three different hyper-heuristics for mine scheduling.~\citet{lamghari2021} combine machine learning (ML), exact algorithms, and heuristic methods to solve the SSOMC and capitalize on the synergies between artificial intelligence and optimization techniques. The authors show that using neural networks as a surrogate model for the downstream subproblem (instead of solving it to optimality) significantly reduces the computational time.

\subsection{Machine learning for combinatorial optimization}

Machine learning has shown promise for tackling small-scale operations research (OR) tasks, such as using graph convolutional neural networks to learn branch-and-bound variable selection policies~\citep{gasse2019}, using hybrid architectures for efficient branching on CPU machines~\citep{gupta2020}.
Moreover, ML has also been used to guide OR algorithms on large-scale tasks, such as using graph neural networks~\citep{morabit2021}, neural networks~\citep{Yaakoubi2019, Yaakoubi2020, morabit2022}, and structured convolutional kernel networks~\citep{Yaakoubi2021} to solve the airline crew scheduling problem~\citep{yaakoubi2019_thesis} as well as the vehicle and crew scheduling problem in public transit and the vehicle routing problem with time windows~\citep{morabit2021,morabit2022}.
When one has a structural understanding of the scheduling problem (i.e., a theoretical understanding of the decisions of the optimization algorithm), ML models can be used to provide fast approximations of these decisions and thus reduce the required computational time. This is referred to in the literature as "learning by imitation": samples of the expected behavior are available and used as demonstrations for learning in the ML model. As such, the ML model has the task of minimizing a cost function that compares the expert's decisions to its own. The application of ML, in this case, can be straightforward, but the assumption that enough demonstration samples are available for learning is a strong assumption that is not always true.

Recent research goes a step further by proposing data-driven methods that can guide OR algorithms to solve combinatorial optimization problems. To solve these problems, state-of-the-art algorithms use various heuristics to make decisions that would otherwise be computationally expensive or difficult to incorporate into the model. Therefore, ML emerges as a fitting predictive tool to make such heuristic scheduling decisions in a more principled and optimized manner~\citep{bengio2020}. As such, following~\citet{katz2017},~\citet{lipets2021} use expert knowledge (also known as domain-knowledge information) to generate a transition probability matrix~\citep{zadorojniy2016} used to define a Markov decision process (MDP). The MDP was used as an orchestrator for scheduling heuristics to solve a complex airline crew scheduling problem (crew pairing and rostering). Concurrently,~\citet{zarpellon2020} present a novel imitation learning framework and hypothesize that the state parameterization of the Branch-and-Bound search tree can help solve mixed-integer linear programming (MILP) problems.~\citet{chmiela2021} propose the first data-driven framework for heuristic scheduling in an exact mixed-integer programming (MIP) solver, where they obtain problem-specific heuristic scheduling that collectively finds many solutions at minimum cost by learning from data describing the performance of primal heuristics.

This line of work is followed by proposing the first data-driven framework for heuristic scheduling in a hyper-heuristic solver that is fully self-managed. Since the SSOMC has no pre-established structural knowledge, ML models can be deployed to discover the problem structure using RL. The model can perform a set of actions in each state, with each action having a reward; the RL model then tries to maximize the expected sum of the rewards. A common challenge in implementing this approach is to find an appropriate reward function so that the learning algorithm does not fall into local minima or terminate without sufficiently exploring the search space. In this manuscript, the use of hyper-heuristics is justified by the extremely large size and complex constraints of the SSOMC and the demonstrated success of simulated annealing~\citep{simulated1,simulated2} and hyper-heuristics~\citep{godoy2003, kumral2013} in solving the SSOMC. Indeed, following the work of~\citet{goodfellow2016, goodfellow2017} and~\citet{lamghari2021}, a self-learning hyper-heuristic called L2P (learn-to-perturb) is proposed.

\section{Problem statement}\label{sec:math_form}

This section presents an overview of the mathematical formulation~\citep{goodfellow2016, goodfellow2017, goodfellowThesis} and the solution approach~\citep{goodfellow2016, lamghari2018} that are used to flexibly model varying mining complexes and solve the SSOMC.

A mining complex is a mineral value chain, where materials extracted from the ground (multiple mines) flow through different processing streams and encounter multiple transformations into final sellable products that are delivered to the mineral market. The goal of the optimization of mining complexes is to provide the market with products, where the supply of materials extracted from the mines is characterized by the uncertainty and variability of the geological attributes. This uncertainty is incorporated into the framework as a group of simulated orebody models with their respective relevant attributes, in the form of geostatistical simulations that generate multiple equiprobable scenarios of the attributes of interest (e.g., grades, densities, material types, etc.)

\begin{figure}
\centering
\includegraphics[width=0.7\linewidth]{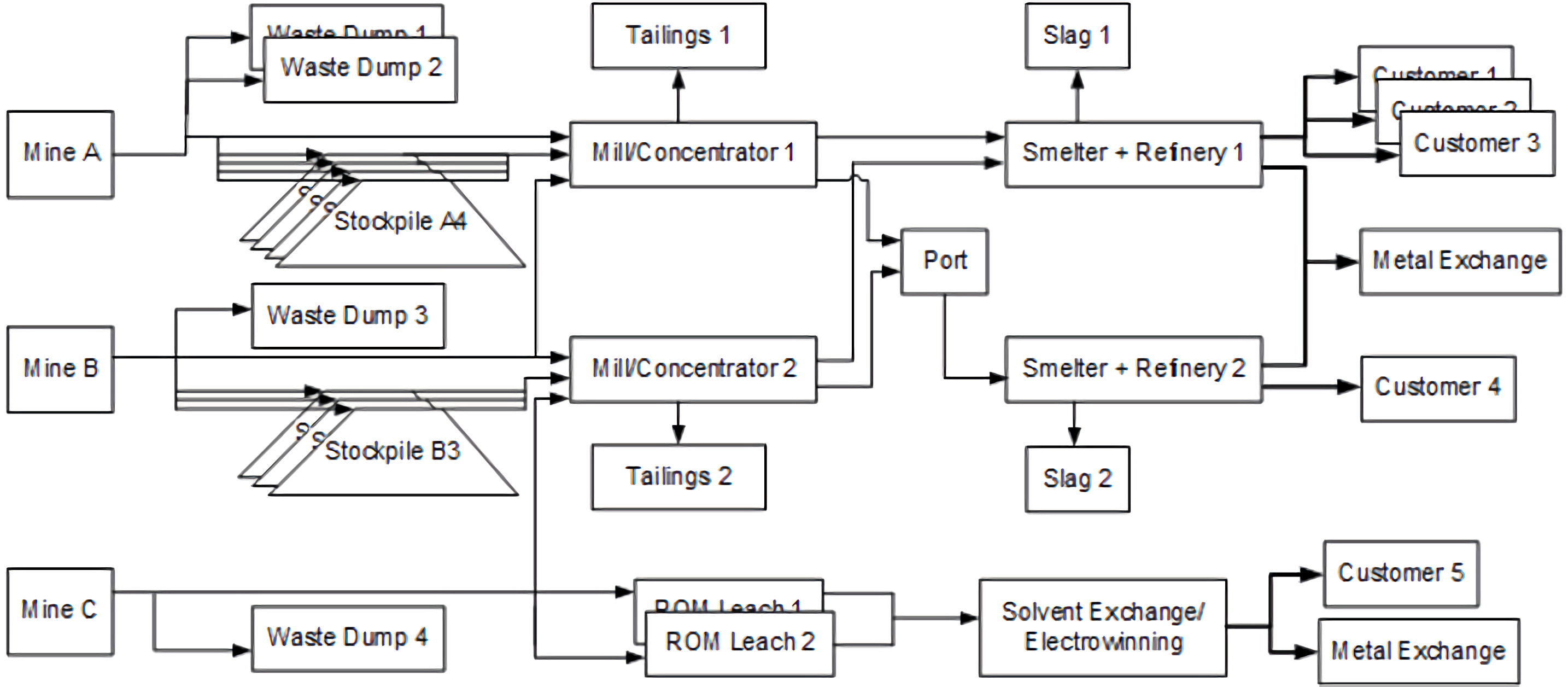}
\caption{An illustration of a mining complex.}
\label{fig:3}
\end{figure}

As in Figure~\ref{fig:3}, the material extracted from mines flows through many different processes in the mining complex, propagating the uncertainty to all operations in the mining complex, which culminates in financial risk. In particular, a mining complex typically contains a set of mines (orebodies), stockpiles (to store incoming materials and meet blending constraints), and processors (to transform incoming products). As such, the SSOMC aims to optimize extraction sequences, destination policies, processing stream decisions, transportation alternatives altogether in one single mathematical model. This framework can be seen as a series of transfer functions that account for the transformations of in-situ material throughout the value chain, which can include various non-linear interactions (blending, non-linear recovery curves, stockpiling, etc.)

In order to optimize the entire mining complex, the problem can be modeled using a two-stage SIP. The objective function maximizes the expected NPV of the mining complex (total expected discounted profit generated from processing the ore) and minimizes the expected recourse costs incurred when the stochastic (scenario-dependent) constraints are violated. Scenario-independent constraints include reserve constraints (ensuring that each block is mined at most once during the horizon), slope and precedence constraints (ensuring that each block can only be mined after all its predecessor blocks have been mined), and mining constraints (imposing a maximum amount of material to be mined per period). Scenario-independent constraints typically include processing constraints (imposing a minimum and maximum amount of material to be processed by each processor per period) and stockpiling constraints (balancing the flow at each stockpile and imposing a maximum amount of ore material to be stockpiled for each stockpile per period).

Due to the exceptionally large size that the model in question can reach and allow for the modeling of varying mining complexes, a generic modeling and solution approach was proposed by~\citet{goodfellow2016, goodfellow2017}. To account for the flow of materials through the mining complex, block attributes are considered to be either primary or hereditary. Primary attributes are those where variables are sent from one location to another (e.g., metal tonnage). Note that these variables can be added directly. Hereditary attributes are variables that are not passed between locations in the mining complex (e.g., costs). These cost functions are assumed to be linear in order to simplify the solution process and to compare the results to the optimal solution (of the relaxed problem) when possible. In the general case, (non)linear formulations can be assigned to these hereditary attributes and are dynamically evaluated during the optimization process.

The flow of materials (and their respective attributes) through the mining complex is defined by three sets of decision variables (which the simultaneous stochastic optimizer can modify), namely mine production scheduling decisions, destination policy decisions, and processing stream decisions. Production scheduling decisions define whether or not a certain mining block is extracted from a given mine in a certain production period. These variables are scenario-independent.

Once mining blocks are extracted, it is necessary to decide where to send each block. As a preprocessing step, multivariate distribution clustering divides blocks into groups based on multiple elements (e.g., grades) in each material type. K-means~\cite{macqueen1967,lloyd1982} is used to define the membership of a certain simulated orebody block for a certain scenario to a given cluster. These membership variables form a robust destination policy that decides where to send all blocks in a cluster (with similar attributes), rather than deciding the destination of individual blocks. Note that a block can be sent to different destinations in various scenarios, but the destination policy decisions are scenario-independent. The optimizer is responsible for deciding where each group (cluster) is sent per period. The use of these multivariate clusters makes the destination decisions more adept at creating destination policies for multi-element mining complexes.

Finally, the flow of materials through the rest of the mining complex is defined using processing stream decisions that define the proportion of a product sent from a location to a destination. These scenario-dependent variables act as recourse actions.
For a detailed formulation of the model, see~\ref{app:formulation}.
In short, it is a two-stage SIP, where first-stage decisions are mine production schedule and destination policies, and second-stage decisions are the processing stream decisions (and penalties for deviations from production targets). The objective function maximizes discounted cash flow and minimizes discounted risk-based deviation penalties, subject to the following constraints: (1) mine reserve and slope constraints, (2) destination policy constraints (clusters of materials are sent to a single destination), (3) processing stream constraints (mass balance for each processor and stockpile, per period, per scenario), (4) attribute calculation constraints (calculate the state of hereditary attributes and quantities in mines), (5) hereditary attribute constraints (tracking deviations from upper and lower bound capacities), (6) stockpile constraints (maintaining quantity balance for each attribute per stockpile per period.)

\section{Solution methodology}\label{sec:method}

This section presents the learn-to-perturb (L2P) solution approach that employs simulated annealing, Tabu search~\citep{glover1998}, and RL to solve the SSOMC. The proposed self-learning hyper-heuristic is constructed by hybridizing a multi-neighborhood simulated annealing algorithm with RL. A probabilistic approach aims at identifying a global optimum by identifying the neighboring states with lower energy or closer to an optimal solution in the solution space.

\subsection{General structure of the learn-to-perturb (L2P) hyper-heuristic}

\begin{figure}
\centering

\resizebox{0.85\linewidth}{!}{
\scalebox{0.7}{
\begin{tikzpicture}   

  \node(poli)[twoblock]{\nodepart[text height=1cm]{one} Policy $\pi$: evaluation \&\ improvement
               \nodepart{two}
                
                \begin{tikzpicture}
                \node(poli6)[crect] {Take action\\$a_t$ = $\pi_{\theta}(s_t)$};
                  \end{tikzpicture}
            
\begin{tikzpicture} \node(poli5)[crect,right= 1mm of poli6] {Update agent's\\weights $\theta$};   
   \end{tikzpicture}};
   
   \node (c5)[circmarker] at (poli.north west) {5};

   \node(past)[twoblock, right=11mm of poli]{\nodepart[text height=1cm]{one} \quad Past Experience
   \nodepart{two} Store experience \\in replay buffer};
   \node (c4)[circmarker] at (past.north west) {4};

 \draw[->,very thick] (past.west) -- (poli.east); 
 
 \node(Reinempty)[below=7mm of poli]{};
 
  \node[draw,fit=(Reinempty) (poli) (past)](Rein) {};
 
  \node(Reintit)[above=0cm of Rein.south, align=center,draw,fill=black!20,inner sep=3mm]{ Reinforcement Learning algorithm };
 %
 %
  \node(samp)[twoblock, below=28.5mm of poli]{\nodepart[text height=1cm]{one} \rule{0mm}{4mm} Sampling function
               \nodepart{two} 
               Update the sampling function based\\
               on past performance and action $a_t$};
 \node (c7)[circmarker] at (samp.north west) {6};             
 

\draw[->,very thick](poli.west) |-  ++(-5mm,0) |- (samp.text west);

  \node(stoch)[twoblock, below= 5mm of samp]{\nodepart[text height=1cm]{one} Stochastic Heuristic: Selection Mechanism
               \nodepart{two}
                
                \begin{tikzpicture}
                \node(stoch1)[crect] {Select a\\heuristic $h_i$\\ };
                  \node (c1)[circmarker] at (stoch1.north west) {1};
                  \end{tikzpicture}
                  
\begin{tikzpicture} \node(stoch2)[crect,right= 1mm of stoch1] {Apply the\\selected heuristic\\ $h_i(x)$};   
\node (c2)[circmarker] at (stoch2.north east) {2};            
   \end{tikzpicture}          };

  \node(simu)[twoblock, right=6mm of stoch]{\nodepart{one} \quad Simulated Annealing Acceptance Criterion
   \nodepart[align=left]{two} \framebox{$x' \gets h_i(x)$} accept the selected heuristic\\[1.2mm]
    \hspace{29mm} OR\\[1.2mm]
     \framebox{$x' \gets x$} \quad reject (restore)};
  \node (c3)[circmarker] at (simu.north west) {3};

   \node (c4)[circmarker] at (past.north west) {4};
    
\draw[->,very thick] (stoch.east) -- (simu.west);

\node(stop)[decision,left= 18 mm of stoch]{Stop?};

\node(yesnode1)[above=10mm of stop.west]{};
\node(yesnode2)[left=-4mm of yesnode1.west]{Yes};

\node(weststop)[left= 9mm of stop]{};

\node(St)[circstate,left= 7 mm of stoch, label=below: {$x$}] {};

\draw[->,dotted, very thick](stop.east)-- node[above]{No}(St.west);
\draw[->,very thick](St.east)--(stoch.west);

\draw[->,very thick] (samp.two west) |-  ++(-15mm,0mm) -|(stop.north);

\draw[->,very thick, dotted](stop.north west)-|++(-2mm,0mm)|- ++(0mm,25mm) node[above](OBs){\begin{tabular}{c}Output\\ Best \\solution
\end{tabular}};

\node(St1)[circstate,right= 5 mm of simu, label=below: {$x'$}] {};
\draw[->,very thick, dotted](simu.east)--(St1);
\draw[->,very thick](St1)|-(past.east);

\node(coll)[draw,very thick, below=21mm of simu.west]{\begin{tabular}{l}Collect problem independent information for domain barrier (e.g., the number of \\heuristics, the changes in evaluation function, a new solution or not, the distance \\between two solutions, etc.)
\end{tabular}};

\node(darr1)[double arrow, inner sep=1mm,  shape border rotate=90, double arrow head extend=1mm, minimum height=6.7mm,draw, below=0mm of stoch.south]{};

\node(darr2)[double arrow, inner sep=1mm,  shape border rotate=90, double arrow head extend=1mm, minimum height=7mm,draw, below=0mm of simu.south]{};

\node(Simulated)[below=7mm of coll]{};

\node[draw,fit= (samp)  (weststop)(simu) (coll) (Simulated) (St1) (OBs)](Simu)  {};
 
 \node(Simutit)[above=0cm of Simu.south, align=center,draw, fill=black!20!white,inner sep=3mm]{ Simulated Annealing };

\node(Doma)[below=3mm of Simu, draw, fill=black!20!white, text width=20cm, text height=2.5mm, align=center, inner sep=3mm]{Domain Barrier};

\node(H1)[circH, below=70mm of St]{$H_1$};
\node(Hn)[circH, right=12mm of H1]{$H_n$};
\node(H2)[circH, below=2mm of H1, xshift=2mm]{$H_2$};
\node(itd)[right=7mm of H2]{\ldots};

\begin{scope}
 \tikzset{every fit/.style={text width=2cm, inner sep=-1mm}}

\node[ellipse,minimum height=2mm, minimum width=51mm,draw,thick, fit=(H1) (H2) (itd) (Hn),xshift=2mm] (HH){};
\end{scope}

\node(PEI)[draw,right=58mm of HH,align=left, thick, text width=50mm]{\\
\hspace{4mm}-- Problem representation\\
\hspace{4mm}-- Evaluation function\\
\hspace{4mm}-- Initial solution ($S_0$)\\
};

\node(S0)[circstate,left= 15 mm of HH, label=below: {$S_0$}] {};

\draw[->,thick] (S0) |- (stop);

\node(Problem)[below=5mm of HH]{};

\node[draw, fit=   (HH) (PEI) (Problem) ](Prob)  {};

 \node(Probtit)[above=0cm of Prob.south, align=center,draw, fill=black!20!white, inner sep=3mm]{Problem Domain};
 

 \end{tikzpicture}
} 
}

   
\caption{Illustration of the learn-to-perturb (L2P) hyper-heuristic, modified from~\citep{burke2013hyper,burke2003hyper,bai2012}.}
\label{Flowchart}
\end{figure}
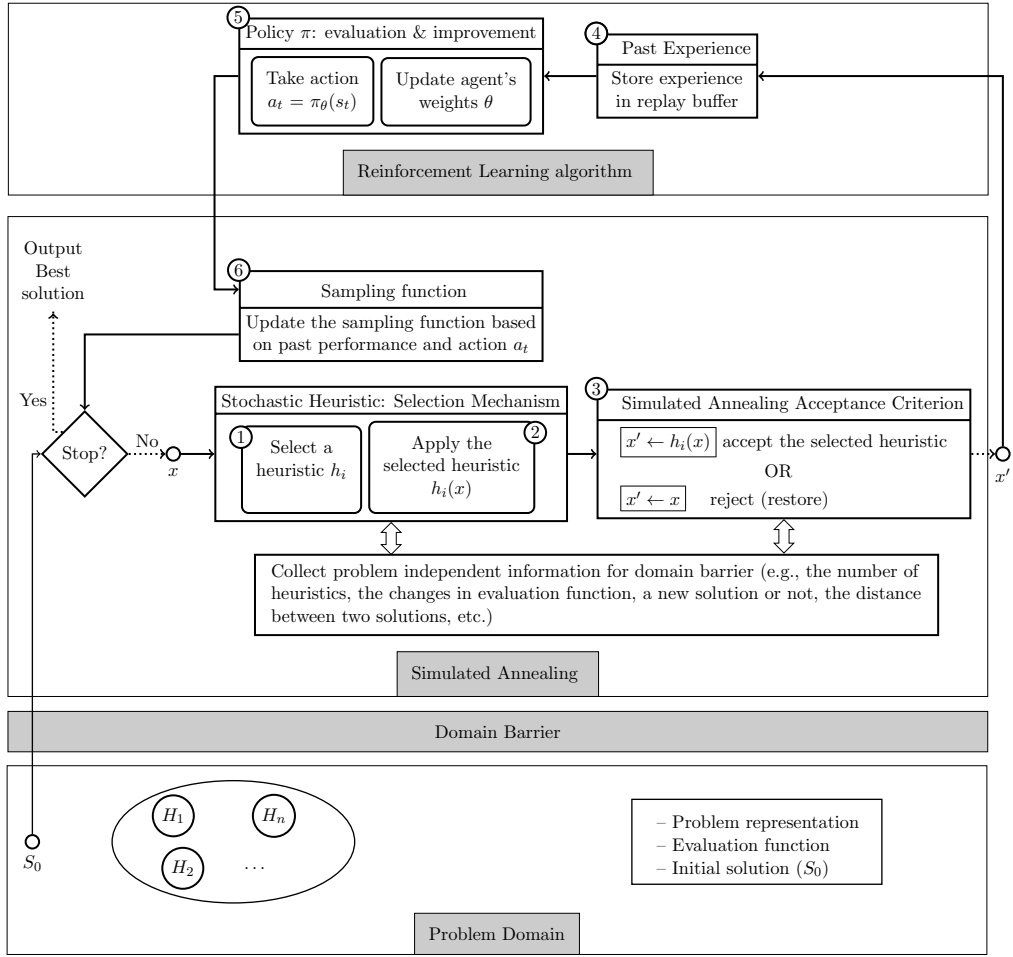

Rather than directly solving an optimization problem, the proposed L2P hyper-heuristic attempts to identify the best heuristic for a given problem state and compensate for the weaknesses of different heuristic methods through simultaneous optimization~\citep{burke2013hyper}. Indeed, while a set of low-level heuristics (perturbations) work on finding optimality simultaneously, a high-level heuristic schedules their sequencing (ordering), ensuring that the best solution is found by combining the perturbations. The proposed L2P hyper-heuristic is illustrated in Figure~\ref{Flowchart} (modified from~\citep{burke2013hyper,burke2003hyper,bai2012}) and is divided into three components that are used simultaneously (rather than in sequential order):
\begin{itemize}
\item \textbf{Problem domain:} includes the problem representation, the evaluation function, the initialization procedure, and the set of low-level heuristics. This is the lower block of Figure~\ref{Flowchart}, and is based on~\citet{goodfellow2017} and~\citet{goodfellowThesis}.
\item \textbf{Simulated annealing algorithm:} consists of a stochastic heuristic selection mechanism, an acceptance criterion and a sampling function, modified from~\citet{lamghari2018}, ~\citet{goodfellow2017}, and~\citet{goodfellowThesis}
\item \textbf{RL algorithm:} predicts the future performance of low-level heuristics to guide the search towards better solutions. This self-managed component learns from past behavior to provide data-driven guidance for searching over the solution space landscape.
\end{itemize}

Following the work of~\citet{lamghari2018},~\citet{goodfellow2017}, and~\citet{goodfellowThesis}, the pseudocode for the L2P hyper-heuristic is presented in Algorithm~\ref{alg:pseudo}. The search is performed using a single candidate solution~\citep{lamghari2018}. At each iteration, L2P improves the current solution $x$ by iterating between heuristic selection and move acceptance (see boxes 1-3 of Figure~\ref{Flowchart}). More precisely, once a heuristic is selected (see box 1), a candidate solution $x$ is modified into a new solution using said heuristic (see box 2). The heuristic is chosen based on its past performance~\citep{lamghari2018} and the expected performance of the RL agent, and simulated annealing~\citep{simulated1,simulated2,simulated3} handles the heuristic selection mechanism and the acceptance criterion.

The algorithm starts with one solution, then moves to a neighboring solution in search of a better solution. This process is repeated until no significant improvement is possible. Heuristic selection is biased toward heuristics that have performed well in previous iterations, and it learns to direct solution sampling strategies to more promising regions of the heuristic space. Specifically, each time a heuristic $h_i$ ($i = 1$, \ldots, $n$, where $n$ is the number of low-level heuristics) is used, the algorithm returns the resulting change in the objective function value $\Delta f(h_i)$ and time required $T(h_i)$. A heuristic is given more importance if it improves solutions quickly per unit time and it is considered less important if it deteriorates the solution and is more computationally expensive. Then, heuristic selection is performed based on the tabu status and performance of the heuristics.

To evaluate the tabu status, if the heuristic does not improve the current solution, it becomes tabu and is temporarily removed from the selection set for $\gamma_T$ iterations where $\gamma_T$ is randomly generated in the interval $[\Gamma_{\min},\Gamma_{\max}]$, and the Tabu list is emptied when all heuristics are made tabu. This strategy allows for better exploration of the heuristic search space because it significantly reduces the probability of repeatedly choosing the same heuristics while leaving the possibility of not banning non-improving heuristics for too long by tuning $\Gamma_{\min}$ and $\Gamma_{\max}$~\citep{lamghari2012,lamghari2018}.

To evaluate the performance of (non-tabu) heuristics, L2P sampels towards the heuristics that are most likely to improve the objective function, by assigning a selection probability $p_i$ computed by normalizing the score function $p_i = \frac{S_F(h_i)}{\sum_{k:H_k}{S_F(H_k)}}$. In what follows, the computation of the score function $S_F$ is explained. Every $\zeta$ iterations, two measures $\pi_1(h_i)$ and $\pi_2(h_i)$ are updated for each heuristic  $h_i$ using Eq.~\eqref{Eq:G2} and Eq.~\eqref{Eq:G3}, respectively~\citep{lamghari2018}. Then, $S_1$ is updated using Eq.~\eqref{Eq:G4}, where $\eta(h_i)$ is the number of times heuristic $h_i$ has been selected in the last $\zeta$ iterations, and $\alpha$ and $\beta$ are two weight adjustment parameters in $[0,1]$ defining the importance given to recent performance and to improving heuristics, respectively. The value of $\beta$ is initialized to 0.5 and is changed depending on whether a new best solution has been found in the last $\zeta$ iterations.

\begin{equation}
\pi_1(h_i) =
\begin{cases}
\pi_1(h_i) + \dfrac{\Delta f(h_i)}{T(h_i)} & \text{if $\Delta f(h_i)$ $>$ $0$ }\\
\pi_1(h_i) & \text{otherwise}
\end{cases}
\label{Eq:G2}
\end{equation}
\begin{equation}
\pi_2(h_i) =
\begin{cases}
\pi_2(h_i) + \dfrac{1}{\Delta f(h_i) T(h_i)} & \text{if $\Delta f(h_i)$ $<$ $0$ }\\
\pi_2(h_i) & \text{otherwise}
\end{cases}
\label{Eq:G3}
\end{equation}
\begin{equation}
S_1(h_i) =
\begin{cases}
S_F(h_i) & \text{if $\eta(h_i) = 0$ }\\
(1-\alpha) S_F(h_i) + \alpha \dfrac{\beta \pi_1(h_i) + (1-\beta) \pi_2(h_i)}{\eta(h_i)} & \text{otherwise}
\end{cases}
\label{Eq:G4}
\end{equation}

\begin{algorithm}
\begin{algorithmic}[1]
\Require
\Statex $x_0$ (initial solution), heuristics $h_i$ ($i = 1$, \ldots , $n$), $k$, $\iter^{\cooling}$ , $\Temp_0$
\Ensure $x_{\best}$
\State Set $x \gets x_0 $, $x_{\best} \gets x$, $\iter \gets 0 $, $\iter' \gets 0 $, and $\Temp \gets \Temp_0$

\textbf{Stage 1: Giving a chance to all heuristics}
\State Add all heuristics to a list $L_H$ (list of not yet selected heuristics)
\While{number of heuristics in $L_H$ $> 0$}
\State Choose a heuristic $h_i$ randomly from the list and generate a new solution $x'$ from $x$ using $h_i$
\If {$x'$ is better than $x_{\best}$}
$x_{\best} \gets x'$
\EndIf
\State Initialize the score function $S_1$
\State Remove $hj$ from the list $L_H$, and set $x \gets x'$
\EndWhile

\textbf{Stage 2: Selecting the heuristics based on their scores and tabu status}
\State Initialize $\pi_1$, $\pi_2$, and $\eta$, and add all heuristics to a list $L_H$ (list of heuristics that are non-tabu)
\While{stopping criterion not met}
\State Select ($h_i$) from $L_H$ based on the normalized probability vector $S_F$ \algorithmiccomment{box 1}
\State Set $\eta(h_i) \gets \eta(h_i) + 1 $
\State Generate a candidate solution $x'$ from $x$ using ($h_i$) \algorithmiccomment{box 2}
\State Evaluate $\Delta f(h_i)$ and $T(h_i)$
\If {$x'$ is better than $x_{\best}$}
Set $x_{\best} \gets x'$ and $\newBest \gets \true$
\EndIf
\If {$\Delta f(h_i) > 0$} $\pi_1(h_i) \gets \pi_1(h_i) + \dfrac{\Delta f(h_i)}{T(h_i)} $
\Else
\If {$\Delta f(h_i) < 0$} $\pi_2(h_i) \gets \pi_2(h_i) + \dfrac{1}{|\Delta f(h_i)| T(h_i)} $
\EndIf
\State Generate a random number $\gamma_T$ in $[ \Gamma_{\min}, \Gamma_{\max}]$ and make $h_i$ tabu for $\gamma_T$ iterations \EndIf
\If {$\Delta f(h_i) > 0$} {$x \gets x'$} \algorithmiccomment{box 3}
\Else \algorithmiccomment{box 3}
\If {$\exp({-\Delta f(h_i)}/{\Temp})> \Uniform (0,1)$}{ $x \gets x'$} \algorithmiccomment{box 3}
\EndIf
\EndIf
\If {$\iter \bmod \iter^{\cooling} = 0$} {$\Temp \gets \Temp * k$} \EndIf
\If {$\iter' < \zeta $}{ $\iter' \gets \iter' + 1$}
\Else
\If {$\newBest = \true$}{ $\beta \gets 1$}
\Else { $\beta \gets \max(\beta - 0.1 , 0)$}
\EndIf
\State Normalize score functions $S_1$, $\pi_1$, and $\pi_2$
\For{$j = 1$, \ldots, $n$}
\State Update the score function $S_1$ using Eq.~\eqref{Eq:G4} and store past experience \algorithmiccomment{box 4}
\State Get new action $S_2$ from RL agent based on past experience ($\pi_1$, $\pi_2$, and $S_1$).  \algorithmiccomment{box 5}
\State Update the score function $S_F$ using Eq.~\eqref{Eq:G5} and initialize $\pi_1$, $\pi_2$, and $\eta$ \algorithmiccomment{box 6}
\EndFor
\State Normalize score function $S_F$, revoke tabu status of all heuristics (reset $L_H$) \algorithmiccomment{box 6}
\State Set $\iter' \gets 1$ and $\newBest \gets \false$
\EndIf
\EndWhile
\State \textbf{return} $x_{\best}$
\end{algorithmic}
\caption{The L2P hyper-heuristic for the simultaneous stochastic optimization of mining complexes, modified from~\citet{lamghari2018} and~\citet{goodfellow2017} and~\citet{goodfellowThesis}.}
\label{alg:pseudo}
\end{algorithm}

The computed score function $S_1$ is used along with $S_2$ (a score function provided by the RL agent) to calculate the final score function $S_F$. To compute $S_2$, all heuristic score measures for the last $L_w$ rounds are provided to the RL agent as input $s_t$ (in the form of a $3D$ matrix). As in box 5 of Figure~\ref{Flowchart}, the agent outputs the action $a_t = \pi_{\theta}(s_t)$, where $\pi_{\theta}$ is the agent's policy that is parameterized by the weights $\theta$. The score function $S_2$ is obtained by normalizing $a_t$. In order to guide the L2P hyper-heuristic to identify the best solution, the RL agent takes a particular action ($a_t$) in each state and transitions from one state to another, and based on the gain or loss occurring, the RL algorithm updates its policy weights in the direction that is likely to maximize the gain (reward)~\citep{bengio2020}. Once $S_1$ and $S_2$ are computed, a voting mechanism (box 6 in Figure~\ref{Flowchart}) is used to compute the final score function $S_F$ based on $S_1$ and $S_2$, as in Eq.~\eqref{Eq:G5}, where $\lambda_{\RL}$ is the RL agent's contribution to the final score function. While it is conceptually possible not to use voting, it is essential for the agent to learn good policy and find good solutions. The voting mechanism combined with sampling typically promotes exploration in RL to ensure the most viable outcome. Note that the L2P variant with $\lambda_{\RL}=0$ will be considered as a baseline in Section~\ref{sec:comp_exp}.

\begin{equation}
S_F(h_i) =  (1-\lambda_{\RL}) S_1(h_i) + \lambda_{\RL} S_2(h_i)
\label{Eq:G5}
\end{equation}

Using the computed $S_F$ to sample the heuristic to apply, the solution vector is updated. The reward $r$ for the agent performing action $S_2$ is $\Delta f(h_i)$ (the change in the objective function value). The agent's performance is ultimately guided by the policy changes and experience gained by the agent over time. At each iteration $t$, before taking the action, the agent stores past experience in a replay buffer. Specifically, as in box 4, the tuple ($S_{t-1}$, $a_{t-1}$, $r_t$, $S_{t}$) is saved in the buffer, which is used to update the agent's policy. And every few iterations (set by a hyperparameter), the agent updates its weights $\theta$ (see box 5).
Note that as the algorithm progresses, the temperature of the simulated annealing selection mechanism is gradually reduced until only minor changes are accepted. The temperature is controlled by the initial temperature ($\Temp_0$), and a cooling schedule, which is defined by a reduction factor, $k \in [0,1)$, and a number of iterations before the reduction factor is applied (see Algorithm~\ref{alg:pseudo}). Note that the is repeated until a predefined stopping condition is met and only problem-independent information flow is allowed between the problem domain and other layers, as is standard in hyper-heuristics~\citep{burke2013hyper}.
An pseudocode for L2P is given in Algorithm~\ref{alg:pseudo}.

\subsection{Reinforcement learning algorithm}\label{sec:reinforcement-learning}

Reinforcement learning (RL) is a goal-oriented algorithm that learns to achieve complex goals and maximize the objective towards a certain dimension, such as winning a game by accumulating the highest possible score in a certain number of moves. In this manuscript, the potential-based reward shaping~\citep{Ng1999, Devlin2012} is used. The potential is obtained with a model-free RL agent (as opposed to a model-based agent) via policy iteration. Policy iteration is used as a form of optimal control to guide the search in the heuristic space, iterating between policy evaluation and improvement. Basically, the current policy is evaluated for each iteration, and the new cost is used to obtain a new improved control policy. These two steps are repeated until the policy improvement step no longer changes the current policy. Referring to Figure~\ref{Flowchart}, policy iterations and the incorporation of past experience are fundamental elements of the RL algorithm implemented in L2P. The notations that are used in this section are reported in Table~\ref{policy_gradient_Notations}.

\begin{table}
	\caption{Notations used in policy gradient algorithms}\label{policy_gradient_Notations}
	\vspace{2mm}
	\centering
	\begin{tabular}{|p{3cm}|p{12.5cm}|} 
		\hline
		Notation & Description\\
		\hline
		$s \; \in \; \mathcal{S}$ & Set of states \\
		\hline
		$a \;\in \; \mathcal{A}$ & Set of actions\\
		\hline
		$r  \;\in \; \mathcal{R}$ & Set of rewards\\
		\hline
		$s_{t},\; a_{t}, \;r_{t}$ & Notation for state, action, and reward at a given time $t$ \\
		\hline
		$ \gamma $ & Discount factor between 0 and 1, can be described as a penalty imposed\\
		 & due to the uncertainty of future rewards\\
		\hline
		$R_{t}$ & Return at time $ t$, $R_t$ = $\sum_{i=0}^{\infty} \gamma^i r_{t+i+1}$ \\
		\hline
		$\prob(s', r \vert s, a)$ & Probability of transitioning from current state $s$ the next state $s'$ by taking \\
		& the action $a$ and receiving the reward $r$ \\
		\hline
		$\pi_{\theta}(a \vert s)$ & policy $\pi$ parameterized with $\theta$ mapping each state with its optimal action\\
		\hline
		$V(s)$ & State value function \\
		\hline
		$V_{w}(.)$ & State value function with parameter $w$ \\
		\hline
		$V^{\pi}(s)$ & State value function while following the policy $\pi$ for a state $s$; \\
		& $V^\pi (s) = \mathbb{E}_{s\sim \pi} [R_t \vert S_t = s]$ \\
		\hline
		$Q^{\pi}(s,a)$ & Action value function computing the expected return from ($s, a$) \\
		& using policy $\pi$ parameterized with parameter $\theta$. It can also be\\
		& referred to as $Q(s,a)$\\
		\hline
		$A(s,a)$ & Advantage function; $A(s,a) = Q(s, a)-V(s)$ \\
		\hline
		$\rho^\pi(s)$ & Stationary distribution of Markov chain for $\pi_\theta$\\
		\hline
	\end{tabular}
\end{table}

\subsubsection{The Markov Decision Process Formulation}\label{sec:mdp}

A RL agent starts from an initial state and updates its weights in the direction of improving the policy, towards an optimal policy that yields the best possible performance.
The agent's learning and state change mechanism are similar to those of humans in that it can be penalized for a series of incorrect decisions/actions and rewarded for correct decisions/actions. In what follows, the policy is denoted by $\pi$, states by $s$, actions by $a$, the discount factor by $\gamma$, and rewards by $r$.

\begin{figure}
	\centering
	\includegraphics[width=0.4\linewidth]{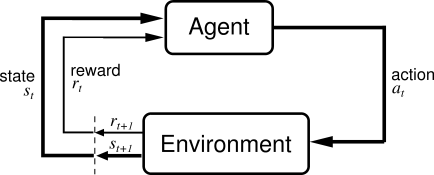}
	\caption{An illustration of the interaction between the RL agent and the environment~\citep{sutton2018}. In the proposed L2P, the agent refers to the reinforcement learning algorithm (box 4-5 in Figure~\ref{Flowchart}).}
	\label{fig:typicalrl}
\end{figure}

As in Figure~\ref{fig:typicalrl}, in a typical RL environment, an agent performs an action $a$ according to policy $\pi$ to transition from the current state $s$ to the next state $s'$. The probability $\prob(s', r \vert s, a)$ is the probability of transitioning from the current state $s_t$ to the next state $s'\sim s_{t+1} $ at time $t$ by taking action $a$ for observation/tuple $(s_t,a_t,r_t,s_{t+1})$. The return $\mathcal{R}_{t}$ at $t$ is: $\mathcal{R}_t = \sum_{i=0}^{\infty} \gamma^i r_{t+i+1}$.

On the one hand, for value-based methods, $\mathcal{R}_t$ is computed to estimate $\pi^*$ that maximizes the expected return $\mathbb{E}_\pi[\mathcal{R}_t]^2$. This is achieved by exploiting value-based solutions, thus computing $V^\pi (s)$ using Eq.~\eqref{eq2}, or by computing $Q^\pi (s,a)$ in Eq.~\eqref{eq3} by exploiting the temporal-difference learning to generate a series of greedy policies by choosing at each time-step the action: $\argmax_{a'}\limits Q_\pi(s,a')$~\citep{srinivasan2018actor}.
\begin{equation}\label{eq2}
    V^\pi (s) = \mathbb{E}_{s\sim \pi} [\mathcal{R}_t \vert \mathcal{S}_t = s]
\end{equation}
\begin{equation}\label{eq3}
    Q^\pi (s,a) = \mathbb{E}_{s\sim \pi} [\mathcal{R}_t \vert \mathcal{S}_t = s, \mathcal{A}_t = a]
\end{equation}

On the other hand, policy gradient methods typically model and update the policy to maximize the rewards. The quality of a policy is measured using the long-term cumulative reward function (formally known as average reward per time-step), as in Eq.~\eqref{eq:1}, where $\rho^\pi(s)$ denotes the stationary distribution of the Markov chain for $\pi_\theta$. Thus, gradient descent is employed on $Y(\pi_\theta)$ to update the weight $\theta$ parameterizing the policy $\pi$.
\begin{equation}\label{eq:1}
	Y(\theta) 
	= \sum_{s \in \mathcal{S}} \rho^\pi(s) V^\pi(s) 
	= \sum_{s \in \mathcal{S}} \rho^\pi(s) \sum_{ a\in \mathcal{A}} \pi_\theta(a \vert s) Q^\pi(s,a)
\end{equation}

\subsubsection{State and Action Encoding}\label{sec:state_encoding}

At each time $t$, the state $s_t$ is defined as all the information (performance) related to all perturbations (low-level heuristics) for the last $L_w$ iterations. Each state is a $3D$ matrix of dimensions $L_w \times 3 \times n$. The first axis (of dimension $L_w$) corresponds to the number of iterations from which useful information is extracted (i.e. the performance of the heuristic). The second axis (of dimension $3$) corresponds to the number of features or measures for each heuristic. For each iteration, the normalized score functions are saved: $\pi_1$ is computed as in Eq.~\eqref{Eq:G2}, $\pi_2$ is computed as in Eq.~\eqref{Eq:G3}, and $S_1$ is computed as in Eq.~\eqref{Eq:G4}. The third axis (of dimension $n$) corresponds to the number of low-level heuristics (see Section~\ref{sec:problem_domain}).

The agent outputs the normalized action ($a_t = \pi_{\theta}(s_t)$), where $\pi$ is the agent's policy parameterized by the weight $\theta$. The output action represents the predicted performance of each heuristic and is guided by updating the parameter $\theta$ of the stochastic policy $\pi_{\theta}$ and the agent's exploration during training. The stochasticity comes from adding a Gaussian noise $N(0,\sigma)$ to the (deterministic) output action $a_t$ to obtain the score function $S_2$. This is used in continuous stochastic control problems (as in this case) to enrich the exploration behaviors. This process can be thought of as adding randomness to the output action in order to explore state and action spaces and thus escape local optima. Although $\sigma$ can be learned during training, it is used as a hyperparameter in this case to avoid an overconfident agent (early convergence of the policy).

After computing the score function $S_2$, a voting mechanism is used to compute the final score function $S_F$ based on $S_1$ and $S_2$, as formulated in Eq.~\eqref{Eq:G5}, where $\lambda_{\RL}$ denotes the contribution of the RL agent to the final score function that will be used for the next optimization iteration. In the next iteration, a heuristic $h_i$ is chosen using the computed score function $S_F$, and a new (updated) solution vector is obtained, along with the resulting change in the value of the objective function $\Delta f(h_i)$ and the time required $T(h_i)$.
Next, the reward given to the agent is computed. Defining a reward is a rather difficult step in designing an RL framework. On the one hand, although defining a sparse reward is easy to encode (e.g., 1 when the problem is solved, 0 otherwise), this setting does not provide a very limited useful learning experience for the agent. On the other hand, when the reward is approximated (e.g., by a linear function approximator instead of solving a downstream process at each iteration), the learning is more likely to fall into local minima~\citep{bengio2020}. In this manuscript, the reward for action $a_t$ is defined as $r_t = \Delta f(h_i)$. Note that to allow generalization for the learned policy across instances, the computed returns are standardized (transformed to have zero mean and unit standard deviation) before using them to update the agent's policy. In addition to $\Delta f(h_i)$, it is also conceptually possible to consider $T(h_i)$ when computing the agent's reward $r_t$. However, this may lead the agent to suggest perturbations that take very little time but provide no improvement in the value of the objective function, resulting in a sub-optimal policy.

\subsubsection{Advantage Actor-Critic (A2C)}

In policy-based methods, the policy is parameterized by $\theta$ (weights for action selection) and $w$ (weights for state evaluation). These parameters are typically updated by gradient descent over $\mathbb{E}(\mathcal{R}_t)$. For example, in the REINFORCE algorithm~\citep{williams1992}, one of the first policy-based methods, the parameters $\theta$ and $w$ are updated in the direction of $\nabla_{\theta}\log \pi(a_{t}|s_t;\theta)$, which can be considered an unbiased estimate of $\nabla_\theta\mathbb{E}(\mathcal{R}_t)$. The main limitation of policy-based methods is that this estimate can have a high variance, which hinders convergence to an optimal policy. To reduce the variance of the estimate while keeping it unbiased, a learned function $b_t(s_t)$ is subtracted from the reward estimate. The resulting equation that defines the gradient is formulated as follows: $\Delta \theta = \nabla_{\theta}\log \pi(a_{t}|s_t;\theta) (\mathcal{R}_t-b_t(s_t))$.

The learned function is referred to as the baseline ($b_t(s_t) \approx V^\pi(s_t)$). In the case where the approximate value function is used as the baseline function, the scaling quantity ($R_t-b_t$) of the gradient can be described as $A (a_t, s_t)=Q(a_t, s_t)-V(s_t)$ where $A (a_t, s_t)$ is referred to as the advantage of an action $a_t$ in a state $s_t$, since $R_t$ is an estimate of $Q(a_t, s_t)$ and $b_t$ is an estimate of $V(s_t)$. The advantage evaluates the quality of an action with respect to the average action for a specific state. This method is called the actor-critic method since the policy $\pi$ acts as an actor and the baseline $b_t$ acts as a critic~\citep{sutton2018, degris2012model}.
The working of A2C is detailed in~\ref{app:rl}. A2C maintains a policy $\pi(a_t,s_t;\theta)$ and the state-value function $V(s_t;\theta_s)$. At each iteration, the algorithm updates the policy and state-value function using n-step returns, for example after reaching a terminal stage or after $t_{max}$ actions. The update can be thought of as $\nabla_{\theta'}\log \pi(a_{t}|s_t;\theta') b(a_{t},s_t;\theta, \theta_s)$ where $b$ can be described as an estimate of an advantage function described in Eq.~\eqref{eq:3}, where the variable $j$ varies from state to state and is bounded by $t_{max}$.
\begin{equation}\label{eq:3}
	b(a_{t},s_t;\theta, \theta_s)
	= \sum_{i=0}^{j-1} \gamma^i r_{t+i}+ \gamma^jV(s_{t+j;\theta_s})-V(s_t;\theta_s)
\end{equation}

The state value function parameter $\theta_s$ and the policy parameter $\theta$ can be shared or completely separate (independent). In our implementation, these two networks share all hidden layers except the last hidden layer and the output layer, in order to speed up learning, while allowing different learning speeds for the actor and value networks. Furthermore, in order to avoid premature convergence, the policy entropy is added to the objective function. The resulting gradient $\Delta_{\theta}$ is expressed in Eq.~\eqref{eq:4}, where $H$ denotes the value of the entropy and $\epsilon$ is a regularization parameter.
\begin{equation}\label{eq:4}
	\Delta_{\theta}	= \nabla_{\theta}\log \pi(a_{t}|s_t;\theta') (\mathcal{R}_t-V(s_t;\theta_s)+\epsilon\nabla_{\theta}, H(\pi(s_t;\theta'))
\end{equation}

\subsubsection{Soft Actor-Critic (SAC)}

The soft actor-critic (SAC)~\citep{sac} is an off-policy actor-critic algorithm that feeds the policy entropy measure $H$ to the reward for robust state exploration while incorporating the maximum possible policy randomness. The training maximizes both $H$ and expected return $ \mathbb{E}_\pi[\mathcal{R}_t]^2$. By using the hyperparameter $\epsilon$ to control the entropy intensity, the reward can be defined as in Eq.~\eqref{eq:5}. Using entropy maximization allows one to choose among optimal and equally good strategies to assign equal probability, and allows further exploration of near-optimal strategies.
%
\begin{equation} \label{eq:5}
	Y(\theta) = \sum_{t=1}^{T}\mathbb{E}_{(s_t,a_{t})\sim \rho_{\pi_\theta}} [r(s_t,a_{t})+\alpha_H H(\pi_{\theta}(.|s_t))]
\end{equation}
\vspace{-5mm}

In practice, SAC learns three functions: (i) the policy $\pi$ with parameter $\theta$ ($\pi_{\theta}$), (ii) the soft Q-value function $Q$ parameterized by $w$ ($Q_{w}$), and (iii) the soft state-value function $V$ parameterized by $\psi$ ($V_\psi$). Using the Bellman equation, the soft Q-value $Q(s_t, a_t)$ and soft state value $V(s_t)$ can be defined as in Eq.~\eqref{eq:6}, where $\rho_{\pi}(s)$ refers to the state and $\rho_{\pi}(s, a)$ refers to the state action pointers of the state distribution provided by the policy $\pi(a|s)$.
\vspace{-2mm}
\begin{equation} \label{eq:6}
	\begin{aligned}
	Q(s_t, a_t) = r(s_t, a_t) + \gamma \mathbb{E}_{s_{t+1} \sim \rho_{\pi}(s)} [V(s_{t+1})], \text{~where~} V(s_t) = \mathbb{E}_{a_t \sim \pi} [Q(s_t, a_t) - \alpha_H \log \pi(a_t \vert s_t)]
	\end{aligned}
\end{equation}
\vspace{-5mm}

First, to update the soft state function (Eq.~\eqref{eq:6}), the MSE (Eq.~\eqref{eq:8}) is minimized, resulting in the gradient in Eq.~\eqref{eq:8_2}, where $\nabla_\psi Y_V(\psi)$ is the gradient and $\mathcal{Z}$ denotes the reply buffer used.
\begin{equation} \label{eq:8}
\begin{aligned}
Y_V(\psi) &= \mathbb{E}_{s_t \sim \mathcal{Z}} \left[\frac{1}{2} \bigl(V_\psi(s_t) - \mathbb{E}[Q_w(s_t, a_t) - \log \pi_\theta(a_t \vert s_t)]\bigr)^2\right]
\end{aligned}
\end{equation}
\vspace{-5mm}
\begin{equation} \label{eq:8_2}
\nabla_\psi Y_V(\psi) = \nabla_\psi V_\psi(s_t)\bigl( V_\psi(s_t) - Q_w(s_t, \theta_t) + \log \pi_\theta (a_t \vert s_t) \bigr)
\end{equation}

Second, to update the soft-Q-function, the Bellman residual is minimized as in Eq.~\eqref{eq:9}, where $\nabla_{w} Y_Q(w)$ is the gradient and $V_{\bar{\psi}}(s_{t+1})$ denotes the target value function. This estimator is updated periodically as an exponential moving average.
\begin{equation} \label{eq:9}
	\begin{aligned}
		Y_Q(w) &= \mathbb{E}_{(s_t, a_t) \sim \mathcal{Z}} [\frac{1}{2}\bigl( Q_{w}(s_t, a_t) - (r(s_t, a_t) + \gamma \mathbb{E}_{s_{t+1} \sim \rho_\pi(s)}[V_{\bar{\psi}}(s_{t+1})]) \bigr)^2] 
	\end{aligned}
\end{equation}
\vspace{-5mm}
\begin{equation} \label{eq:9_2}
		\nabla_{w} Y_Q(w) = \nabla_{w} Q_{w}(s_t, a_t) \bigl( Q_{w}(s_t, a_t) - r(s_t, a_t) - \gamma V_{\bar{\psi}}(s_{t+1})\bigr) 
\end{equation}

Third, SAC also minimizes the Kullback-Leibler divergence (relative entropy) by employing Eq.~\eqref{eq:10}. The update ensures that $Q^{\pi_\text{new}}(s_t, a_t) \geq Q^{\pi_\text{old}}(s_t, a_t)$, where $\Pi$ denotes the set of potential policies that can be used to model our policy. $D^{\pi_\text{old}}(s_t)$ can be thought of as an intractable partition function that is used to normalize the distribution, but it never contributes to the calculation of the gradient. The working of the SAC algorithm is detailed in~\ref{app:rl}.
\begin{equation} \label{eq:10}
	\begin{aligned}
		\pi_\text{new} 
		&= \arg\min_{\pi' \in \Pi} \mathcal{Z}_\text{KL} \left( \pi'(.\vert s_t) \| \exp(Q^{\pi_\text{old}}(s_t, .) - \log D^{\pi_\text{old}}(s_t)) \right) \\[6pt]
		\text{objective for update: } Y_\pi(\theta) &= \nabla_\theta \mathcal{Z}_\text{KL} \left( \pi_\theta(. \vert s_t) \| \exp(Q_{w}(s_t, .) - \log D_{w}(s_t)) \right) \\[6pt]
		&= \mathbb{E}_{s_t\sim\pi} [ \log \pi_\theta(a_t \vert s_t) - Q_{w}(s_t, a_t) + \log D_{w}(s_t) ]
	\end{aligned}
\end{equation}

\subsubsection{Proximal Policy Optimization (PPO)}

The motivation behind incorporating proximal policy optimization into L2P comes from the fact that policy gradient algorithms can be unstable in continuous control problems, as is the case for the current RL environment. In order to improve the training process and stability and to avoid diverging from the optimal policy, the trust region policy optimization (TRPO)~\citet{schulman2015trust}, used the Kullback-Leibler divergence constraint~\citep{hu2013kullback} on the policy so that difference between the updated policy and the previous policy is relatively small.

Since this procedure is computationally too complex,~\citet{ppo} propose proximal policy optimization (PPO), which retains the efficiency of TRPO~\citep{schulman2015trust} while being computationally feasible (tractable). PPO circumvents the Kullback-Leibler divergence constraint by employing a clipped surrogate objective. That is, PPO keeps $r(\theta)$ (the ratio of $\theta_{old}$ to $\theta$, i.e., $\frac{\pi_\theta(a \vert s)}{\pi_{\theta_\text{old}}(a \vert s)}$) within a small interval centered on 1, such as $[1-\epsilon,1+\epsilon]$, with $\epsilon$ being a hyperparameter.

The clipped advantage function becomes as in Eq.~\eqref{eq:13}, where the clipping function $\CLIP(r(\theta),1-\epsilon,1+\epsilon)$ retains the ratio between $(1+\epsilon)$ and $(1-\epsilon)$. Note that the estimated advantage function $\hat{\mathcal{W}}(.)$ is used here instead of true advantage function $A(.)$ because true rewards are unknown in this case. The objective function to be maximized in order to update the weights $\theta$ becomes as in Eq.~\eqref{eq:14}. The objective function employs an entropy term for sufficient state exploration and $p_1, p_2$ are hyper-parameters controlling the exploration/exploitation tradeoff~\citep{ppo}.
\begin{equation} \label{eq:13}
Y^{\CLIP} (\theta) = \mathbb{E} [ \min( r(\theta) \hat{\mathcal{W}}_{\theta_{\old}}(s, a), \CLIP(r(\theta), 1 - \epsilon, 1 + \epsilon) \hat{\mathcal{W}}_{\theta_{\old}}(s, a))]
\end{equation}
\vspace{-4mm}
\begin{equation} \label{eq:14}
	Y^{\CLIP'} (\theta) = \mathbb{E} [ Y^{\CLIP} (\theta) - p_1 (V_\theta(s) - V_\text{target})^2 + p_2 H(s, \pi_\theta(.)) ]
\end{equation}
Following~\citet{hsu2020revisiting}, to avoid instability in continuous action spaces and (over)sensitivity to initialization (e.g., if the optimal local actions are close to initialization), the probability ratio is clipped to regularize the policy and a continuous Gaussian distribution based on the parameterized policy is used as in Section~\ref{sec:state_encoding}. The working of our PPO implementation is shown in~\ref{app:rl}, where the surrogate loss is optimized using the Adamax optimizer~\citep{kingma2014} (or any similar optimizer).

\subsection{Problem domain component}\label{sec:problem_domain}

To produce new solutions, the hyper-heuristic described in the previous section uses 38 simple low-level perturbative heuristics from~\citet{goodfellow2016,goodfellow2017,goodfellowThesis}, referred to as $h_i$ ($i = 1$, \ldots, $n$, where $n$ is the number of low-level heuristics). Each heuristic aims to improve the current solution and examines a subset of one of the following four neighborhoods:

\begin{enumerate}
\item  \textbf{Block extraction sequence perturbations~\citep{goodfellow2016,goodfellow2017}}:
For each of this neighborhood's heuristics, a block is randomly selected. The block's extraction period is changed, according to the heuristic's rules, and any potentially violated slope constraints  are fixed. Specifically, (1) the current block that is visited is set as the central block. (2) If there are slope constraint violations for any of the blocks in the predecessor/successor set of the central block when the central block is moved to the new time period, that block will be added to the set of blocks whose sequences are to be changed. (3 - a) If the previous envelope of step~2 has any violations of the slope constraint, the current block is set as the block that overlies the current block, and step~2 is repeated. (3 - b) The extraction period of the previously identified blocks is changed if there are no slope constraint violations in the envelope of step~2 or if there is no overlying block (air block). Note that a predecessor of block $i$ is a block that must be extracted to have access to $i$ and block $j$ is a successor of block $i$ if and only if $i$ is a predecessor of $j$. If the extraction period of a block $i$ is advanced (resp. delayed), this process is equivalent to the predecessors (resp. successors) of $i$ being moved to the same new extraction period as $i$. The set formed by a block $i$ and its predecessors (resp. successors) mined in the same period is called the inverted cone (resp. cone), whose base (resp. apex) is $i$.

\item \textbf{Cluster destination policy perturbations~\citep{goodfellow2016}:}
A cluster is a group of blocks in a mine that can be grouped together based on certain characteristics. Clusters can be formed by grouping blocks with similar block attributes. The membership of a block in a cluster can change throughout simulations as block characteristics vary in different simulations. The overall optimal solution is essentially the sum of the optimal results from all simulations (scenarios). In the heuristics of this neighborhood, a cluster destination decision variable is randomly selected and sent to a different destination, if possible.

\item \textbf{Destination policy perturbations~\citep{goodfellowThesis}:}
A destination policy is often employed to group blocks into bins of ore and waste blocks. In this way, optimization is limited to looking at a collection of blocks instead of individual blocks. For example, in a time-varying grade policy, each block characteristic is assigned a value below which mining becomes infeasible. These heuristics will dynamically vary the cut-off range characteristics of clusters when cluster-related perturbations occur, allowing for further optimization.

\item \textbf{Processing stream perturbations~\citep{goodfellow2017,goodfellowThesis}:}
Open-pit mines can operate with one or multiple process streams, where materials are extracted in bulk, processed, and sent to specific destinations. To find these destination policies, the process stream variables are optimized to maximize the overall economic value of the operation. Process streams can have material quality constraints, although the current model uses the cut-off grade constraint. The transfer of materials between two locations is typically dependent on the process stream configuration and is, therefore, an integral part of the modeling and optimization process. A random normal number $N(y_{i,j,t,s}, 0.1)$ is used to change a randomly selected process stream variable. Keeping the distribution's variance small enough facilitates both global and local search.
\end{enumerate}

\section{Computational experiments}\label{sec:comp_exp}

\subsection{Experimental setting, benchmark instances and parameter setting}\label{sec:bench_instances}

To evaluate the performance of the proposed L2P hyper-heuristic, numerical experiments on five instances (of the simultaneous stochastic optimization of mining complexes) of various sizes and characteristics are reported, divided into three different stages of testing. The benchmark instances are briefly described below and detailed in Table~\ref{table:1}. Different variants of the L2P hyper-heuristic are compared based on the evolution of the objective function value with respect to the number of iterations and the computational time. L2P--Baseline refers to the L2P hyper-heuristic with $\lambda_{\RL} = 0$ (not using RL). L2P--A2C, L2P--PPO, and L2P--SAC refer to the L2P hyper-heuristic variants using A2C, PPO, and SAC as RL agents, respectively. All variants use $\lambda_{\RL} = 0.5$. The L2P variants are compared based on on the evolution of the objective function value with respect to the number of iterations and the computational time.

In the first stage of testing, three instances (mining complexes) $I_1$-$I_3$ are used. This stage is intended to make an initial evaluation of the proposed L2P hyper-heuristic with emphasis on the generalization capability of the proposed method. The instances are small size with up to 15,000 blocks for one mine or up to 4,000 blocks per mine for two mines. Note that online learning (weight updating) is disabled in this stage, so the resulting performance is only obtained from the pre-training stage. Production, processing, and stockpile capacities are fixed for all periods/years. $I_1$ contains a copper mine (one type of metal), a processor, and a waste dump. $I_2$ contains two mines with two types of copper materials (oxides and sulfides), two processors (oxide leach pad and sulfide processor), and a waste dump. $I_3$ is similar to $I_2$, except that it has $\times$2 times as many blocks and a sulfide stockpile. Each time, the method is initialized, pre-trained on one of the instances, and then tested on all three. $I_3$ is similar to $I_2$, except that it has $\times 2$ times as many blocks and a sulfide stock. To test, each L2P variant is trained using one instance, and then tested using all three instances. Since the L2P starts with a random initial solution and the heuristic selection mechanism is stochastic, all experiments are run 10 times and the results in this stage are reported as mean and standard deviation. To compare the different L2P variants, since finding the optimal value for any of the instances is not possible, the linear relaxation of the model is used, which leads to a weak bound on the objective value of the optimal solution. A time limit of four weeks is set to solve the linear relaxation for each instance using CPLEX 12.10.

\begin{table}
\centering
\caption{Overview of the instances in the five benchmark instances.}\label{table:1}
\vspace{2mm}
\resizebox{\linewidth}{!}{ 
\begin{tabular}{ l | l l l || l l }
\hline
 Instance & $I_1$ & $I_2$ & $I_3$ & $I_4$ & $I5$ \\
\hline
Number of mines & 1 & 2 & 2 & 1 & 2 \\
Number of blocks ($N$) & 14,952 & 2,114 and 1,123 & 3,234 and 3,783 & 213,003 & 136,236 and 108,753 \\
Number of periods ($T$) & 8 & 20 & 20 & 8 & 20 \\
Number of simulations ($S$) & 20 & 20 $\times$ 20 & 20 $\times$ 20 & 20 & 20 $\times$ 20 \\
Number of processors ($P$) & 1 & 2 & 2 & 1 & 2 \\
Number of stockpiles & 0 & 0 & 1 & 0 & 1 \\
Number of waste dumps & 1 & 1 & 1 & 1 & 1 \\
Metal type & Copper & Oxides and Sulfide & Oxides and Sulfide & Copper & Oxides and Sulfides \\
Block weight ($w_i$) in tonnes & 10,800 & 16,500 and 16,560 & 16,500 and 16,560 & 10,800 & 16,500 and 16,560 \\
Upper bound on mining ($\times 1000$) & 400 & 2,000 and 2,000 & 2,000 and 2,500 & 28,000 & 30,000 and 25,000 \\
Lower bound on processing ($\times 1000$) & 0 & 1,400 & 0 and 3,200 & 7,250 & 0 and 34,000 \\
Upper bound on processing ($\times 1000$) & 300 & 1,800 & $\infty$ and 4,000 & 7,500 & $\infty$ and 36,000 \\
Upper bound on stockpile capacity ($\times 1000$) & 0 & 0 & 2,000 & 0 & 8,000 \\
Discount rate for cash flow $d_1$ & 0.1 & 0.1 & 0.1 & 0.1 & 0.1 \\
Discount rate for recourse costs $d_2$ & 0.1 & 0.15 & 0.15 & 0.1 & 0.15 \\
Discount rate for recourse costs $d_3$ & 0.1 & 0.15 & 0.15 & 0.1 & 0.15 \\
Discount rate for recourse costs $d_4$ & 0.1 & 0.15 & 0.15 & 0.1 & 0.15 \\
Discount rate for recourse costs $d_5$ & 0.1 & 0.15 & 0.15 & 0.1 & 0.15 \\
Discount rate for recourse costs $d_6$ & 0 & 0 & 0.15 & 0 & 0.15 \\
unit shortage cost (processors) & 8 & 7 & 7 & 8 & 7 \\
unit surplus cost incurred (processors) & 5 & 10 & 10 & 5 &  10 \\
unit surplus cost incurred (mining)  & 1 & 15 and 6.5 & 15 and 6.5 & 1 & 15 and 6.5 \\
unit surplus cost incurred (stockpiling)  & 0 & 0 & 10 & 0 & 10 \\
\hline
\end{tabular}
}
\end{table}

In the second stage of testing, the agents are trained and tested on the large-scale instance of the SSOMC $I_4$. The instance $I_4$ is similar to $I_1$ except that it has $\times 14$ more blocks and imposes a lower processing limit, and this stage is intended to evaluate the performance of the L2P hyper-heuristic on a realistic case study. Each L2P variant is trained using $I_4$, then tested on the same instance $I_4$. Since CPLEX was unable to solve the linear relaxation of $I_4$ within the time limit (four weeks), the objective value of the best solution found (after several trials) is therefore used as a comparison. This point of reference is referred to as $Z^*$. In this stage, all variants of the L2P hyper-heuristic will be compared based on the evolution of the objective function value (with respect to $Z^*$ in \%) as a function of the number of iterations and the execution time (in minutes), in the form of confidence intervals for the estimated P10, P50, and P90 quantiles.

\begin{table}
\centering
\caption{Parameters used for the L2P hyper-heuristic algorithm.}\label{table:2}
\vspace{2mm}
\resizebox{0.7\linewidth}{!}{ 
\begin{tabular}{ l l }
\hline
 Parameter description & Value \\
\hline
Number of layers in critic network & 3 \\
Number of units in critic network & 200, 200, 1 \\
Optimizer for critic network & Adamax \\
Critic learning rate & 0.0001 \\
Number of layers in actor network & 3 \\
Number of units in actor network & 200, 200, 38 \\
Optimizer for actor network & Adamax \\
Actor learning rate & 0.001 \\
Number of steps before updating the network & 5 \\
Standard deviation for Gaussian noise & 0.05 \\
Window length ($L_w$) & 5 \\
Clip gradient norm & 1 \\
Discount factor & 0.9 \\
Agent's contribution ($\lambda_{\RL})$ & 0.5 \\
Number of iterations for training & $10^6$ ($\approx$ $12$--$16$ hours) \\
\hline
\end{tabular}
}
\end{table}

In the third stage of testing, the large-scale instance $I_5$ is used to evaluate the performance of the L2P hyper-heuristic on a realistic case study. The instance $I_5$ is similar to $I_3$ except that it has $\times 75$ times more blocks. Each L2P variant is pre-trained on the instance $I_4$, then the pre-trained weights are transferred and the variant is tested on the instance $I_5$. Note that for both the second and third stages of testing, online learning is enabled as it would be for an end-user. The only difference with the second stage of testing is that the pre-training (warm-starting the policies) is performed here on a completely different instance.

Table~\ref{table:2} reports parameters of the L2P hyper-heuristic. In particular, the actor and critic networks use three layers containing 200 neurons per hidden layer, and the Adamax optimizer~\citep{kingma2014} is used to update the networks' weights using $(10^{-3})$ and $(10^{-4})$ as learning rates for the actor and critic, respectively. For fast execution and efficient memory and GPU allocation, while the mining complex framework and simulated annealing algorithm are implemented in C++, the RL agents are written in Python using the Pytorch library. The experiments are performed on a standard Windows machine equipped with an 8-core processor, 32 GB of RAM, and a GPU (4 GB of GDDR6 memory).

\subsection{Numerical results}\label{sec:num_res}

\subsubsection{First stage of testing}

Using the value of the obtained linear relaxation, Tables~\ref{nt2}-\ref{nt1} report the execution time ($T_{gap=1\%}$) and the number of iterations ($\Iter_{gap=1\%}$), respectively, required by the variants to obtain a solution whose objective function value is $1\%$ away from the linear relaxation value. The last row of Table~\ref{nt2} reports the time required by CPLEX to solve the linear relaxation. The best results obtained for each instance are indicated in bold.


\begin{table}
\centering
\caption{Result summary (mean and standard deviation) of the number of iterations $\Iter_{\gap =1\%}$ $(\times 1000)$ to achieve a solution that is $1\%$ far from the linear relaxation}\label{nt2}
\vspace{2mm}
\begin{tabular}{|c|c|c|c|c|c|}
\hline
\multicolumn{1}{|c|}{} & \multicolumn{1}{c|}{} & \multicolumn{3}{c}{$\Iter_{\gap =1\%}$ $(\times 1000)$} & \multicolumn{1}{c|}{} \\
\hline
\pbox{20cm}{Instance used \\for testing} & \pbox{20cm}{Instance used \\for training} & \pbox{5cm}{L2P--Baseline} & L2P--A2C & L2P--SAC & L2P--PPO\\
\hline
 & $I_1$ & & \textbf{3.9} (1.1) & 6.5 (0.8) & 4.5 (\textbf{0.6})\\
\cline{2-2}\cline{4-6}
$I_1$ & $I_2$ & 8.4 (\textbf{0.3}) & 4.6 (1.3) & 6.8 (0.8) & 5.1 (0.7)\\
\cline{2-2}\cline{4-6}
& $I_3$ & & 4.4 (1.4) & 6.7 (0.9) & 5.2 (0.7)\\
\hline
& $I_1$ & & 39.4 (8.1) & 64.2 (3.8) & 45.9 (3.6)\\
\cline{2-2}\cline{4-6}
$I_2$ & $I_2$ & 78.4 (\textbf{2.5}) & \textbf{35.2} (7.2) & 62.4 (3.4) & 41.0 (\textbf{3.2})\\
\cline{2-2}\cline{4-6}
& $I_3$ & &36.5 (7.0) & 64.8 (3.5) & 40.7 (3.3)\\
\hline
& $I_1$ & & 69.6 (14.3) & 113.5 (6.7) & 81.1 (6.4)\\
\cline{2-2}\cline{4-6}
$I_3$ & $I_2$ & 138.6 (\textbf{4.4}) & 64.2 (12.7) & 110.3 (6.0) & 72.5 (\textbf{5.7})\\
\cline{2-2}\cline{4-6}
& $I_3$ & & \textbf{62.5} (12.4) & 114.6 (6.2) & 72.0 (5.8)\\
\hline
\end{tabular}
\end{table}

The results show that all RL-based variants of the L2P hyper-heuristic outperform L2P--Baseline, where L2P--Baseline refers to the L2P hyper-heuristic with $\lambda_{\RL} = 0$ (not using RL). On average, L2P--A2C outperforms the other hyper-heuristic variants, reducing the number of iterations and execution time by 45\%-55\%. L2P--PPO reduces the number of iterations and execution time by 38\%-48\%. And L2P--SAC reduces the number of iterations and execution time by 17.3\%-22.6\% and 17\%-22\%, respectively. When comparing the performance of variants trained on different instances, the slightly better performing variant is naturally the one tested on the same instance used for training. Although training is performed on the heuristic space rather than the solution space, the heuristic choice and performance may vary slightly from instance to instance, which is a good indicator of the generalization capacity of the proposed method, especially since the weights are not updated in real-time at this stage of testing (i.e., once the weights are pre-trained and transferred, no online learning is used). This difference in performance may only become smaller in an actual use case scenario, where online learning is enabled.


\begin{table}
\centering
\caption{Result summary (mean and standard deviation) of the computational time $T_{\gap=1\%}$ (in minutes) to achieve a solution that is 1\% far from the linear relaxation (LR)}\label{nt1}
\vspace{2mm}
\begin{tabular}{|c|c|c|c|c|c|c|}
\hline
\multicolumn{1}{|c|}{} & \multicolumn{1}{c|}{} & \multicolumn{3}{c}{$T_{\gap=1\%}$ (in minutes)} & \multicolumn{1}{c|}{} & \multicolumn{1}{c|}{LR} \\
\hline
\pbox{20cm}{Instance used \\for testing} & \pbox{20cm}{Instance used \\for training} & \pbox{5cm}{L2P--Baseline} & L2P--A2C & L2P--SAC & L2P--PPO & \\
\hline
 & $I_1$ & & \textbf{8.5} (2.4) & 14.2 (1.7) & 9.8 (\textbf{1.3}) & \\
\cline{2-2}\cline{4-6}
$I_1$ & $I_2$ & 18.3 (\textbf{1.2}) & 10.0 (2.8) & 14.8 (1.7) & 11.1 (1.5) & 2100 \\
\cline{2-2}\cline{4-6}
& $I_3$ & & 9.6 (3.0) & 14.6 (1.9) & 11.3 (1.5) & \\
\hline
& $I_1$ & & 32.7 (6.7) & 53.4 (3.2) & 38.2 (3.0) & \\
\cline{2-2}\cline{4-6}
$I_2$ & $I_2$ & 65.2 (\textbf{1.9}) & \textbf{29.3} (5.9) & 51.9 (2.8) & 34.1 (\textbf{2.7}) & 6480 \\
\cline{2-2}\cline{4-6}
& $I_3$ & &30.4 (5.8) & 53.9 (2.9) & 33.8 (\textbf{2.7}) & \\
\hline
& $I_1$ & & 59.2 (12.2) & 96.1 (5.7) & 68.9 (5.4) & \\
\cline{2-2}\cline{4-6}
$I_3$ & $I_2$ & 117.8 (\textbf{3.7}) & \textbf{52.9} (10.8) & 93.8 (5.1) & 61.6 (\textbf{4.8}) & 8400 \\
\cline{2-2}\cline{4-6}
& $I_3$ & & 54.8 (10.5) & 97.4 (5.3) & 61.2 (4.9) & \\
\hline
\end{tabular}
\end{table}

Results also show that although L2P--A2C outperforms the other variants, its performance is not as consistent across runs as L2P--PPO (the second-best performing variant). In particular, the standard deviations (both for $\Iter_{gap=1\%}$ and $T_{gap=1\%}$) of L2P--A2C are larger than those of L2P--PPO. This may be explained by the fact that PPO takes much smaller steps to update its policy, and thus presents a more stable strategy than A2C (albeit less efficient). In addition, the standard deviation of L2P--Baseline is smaller than the standard deviation of all other variants. This is because, although the choice of heuristic is the result of a sampling process, the score function $S_1$ is biased towards the best performing heuristic, and $S_2$ might be more biased towards exploration, especially if the state given by the recent performance has not been seen before. Thus, using the score function $S_1$ solely may produce more ``deterministic'' (and less noisy) behavior.

Due to the large size of the instances in question, CPLEX takes much longer (150-250 times) to solve the linear relaxation than any of the L2P hyper-heuristic variants take to solve the instances. Finding integer solutions using CPLEX would be impossible, let alone adding more components to the mining complexes presented in this manuscript. In summary, the best performing variant of the L2P hyper-heuristic is L2P--A2C, both for the instance on which it was trained and for the instances on which it was not trained. The RL--based variant with the smallest standard deviation is L2P--PPO, and the least-performing RL--based variant is L2P--SAC.

\subsubsection{Second stage of testing}

\begin{figure}
\begin{minipage}[t]{.45\textwidth}
\centering
\includegraphics[width=\textwidth]{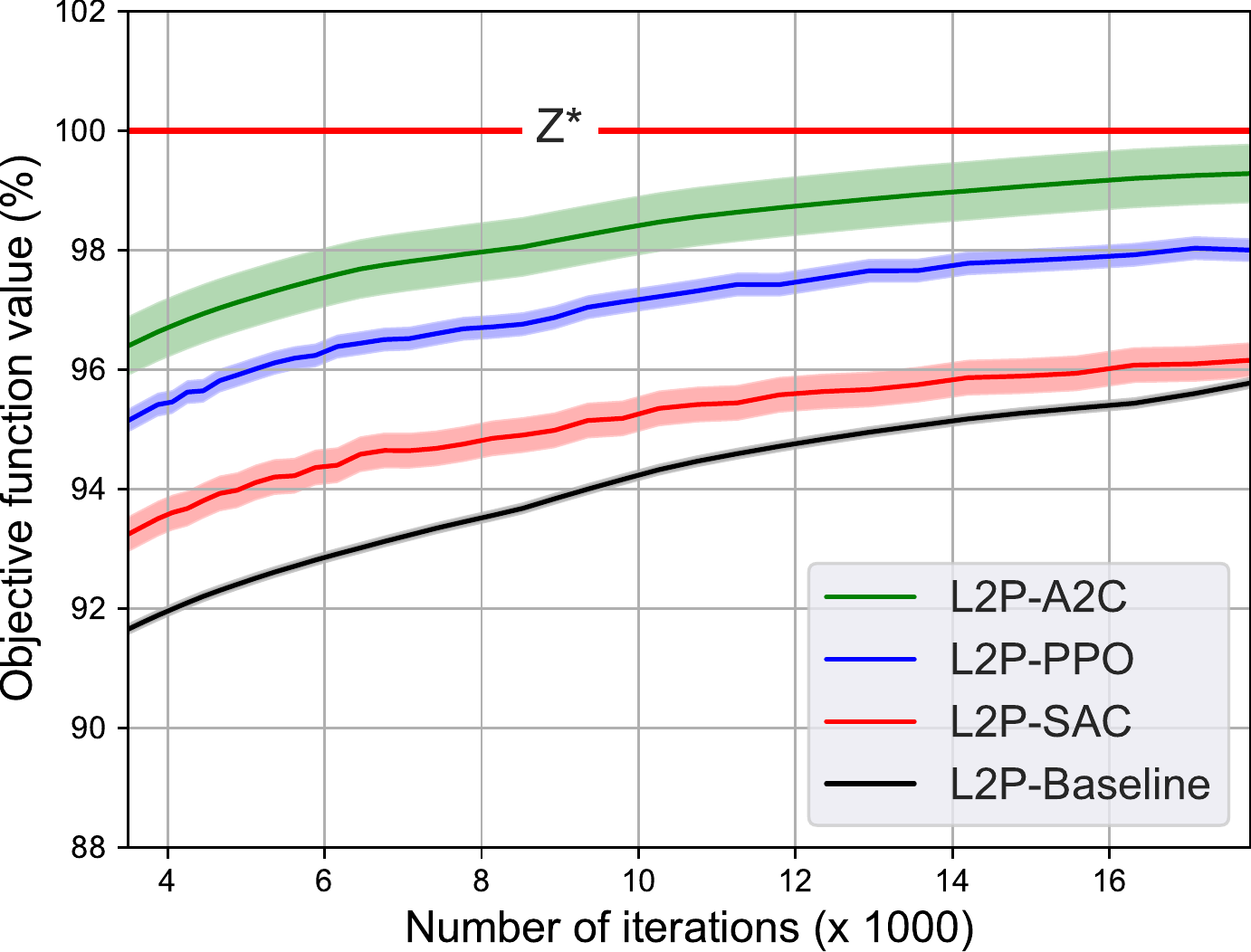}
\end{minipage}
\hfill
\begin{minipage}[t]{.45\textwidth}
\centering
\includegraphics[width=\textwidth]{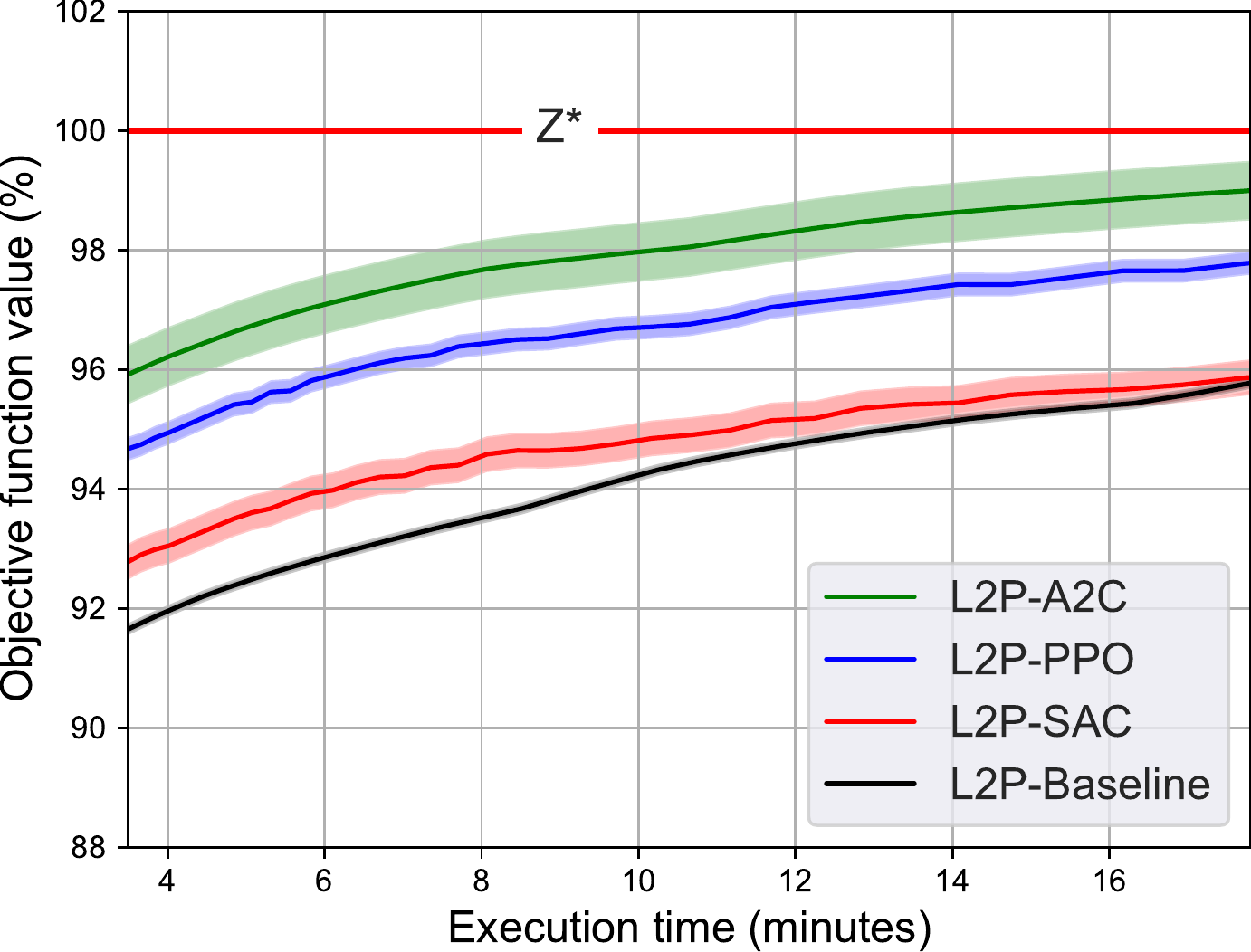}
\end{minipage}
\caption{Evolution of the objective function value (in \%) for Instance $I_4$ with respect to the number of iterations ($\times1000$) (left plot) and the execution time (in minutes) (right plot).
}
\label{fig:day2}
\end{figure}

As shown in Figure~\ref{fig:day2}, L2P--A2C is by far the best performing variant of the L2P hyper-heuristic, reducing the number of iterations by $50$--$80\%$ and the execution time by $40$--$75\%$. The second-best candidate is L2P--PPO, which reduces the number of iterations by $45$--$65\%$ and the execution time by $40$--$65\%$. Compared to the other variants, L2P--PPO has the smallest standard deviation in terms of the number of iterations and execution time. L2P--SAC is the least performing variant, showing that SAC is not well suited for the current environment. This can be explained by two reasons. First, SAC uses entropy regularization to prevent premature convergence. This can make SAC take much longer to learn compared to A2C and PPO. Second, the trade-off between exploration and exploitation within SAC is controlled by the coefficient $\alpha_H$, which needs to be changed from one environment to another and requires careful tuning. This makes SAC difficult to include in an off-the-shelf solver with little to no pre-tuning and pre-training.

\begin{figure}
\begin{minipage}[t]{.4\textwidth}
\centering
\includegraphics[width=\textwidth]{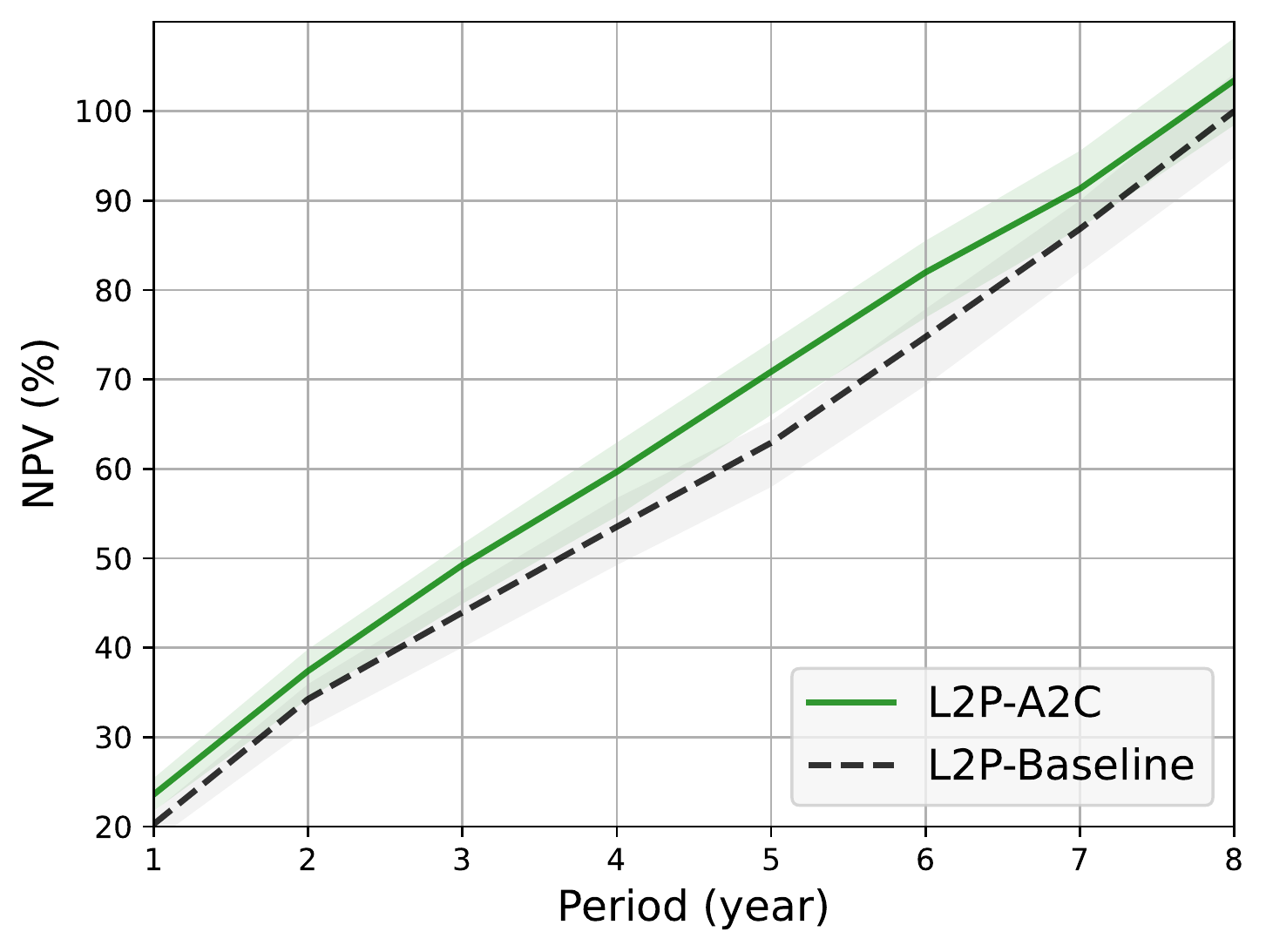}
\end{minipage}
\hfill
\begin{minipage}[t]{.4\textwidth}
\centering
\includegraphics[width=\textwidth]{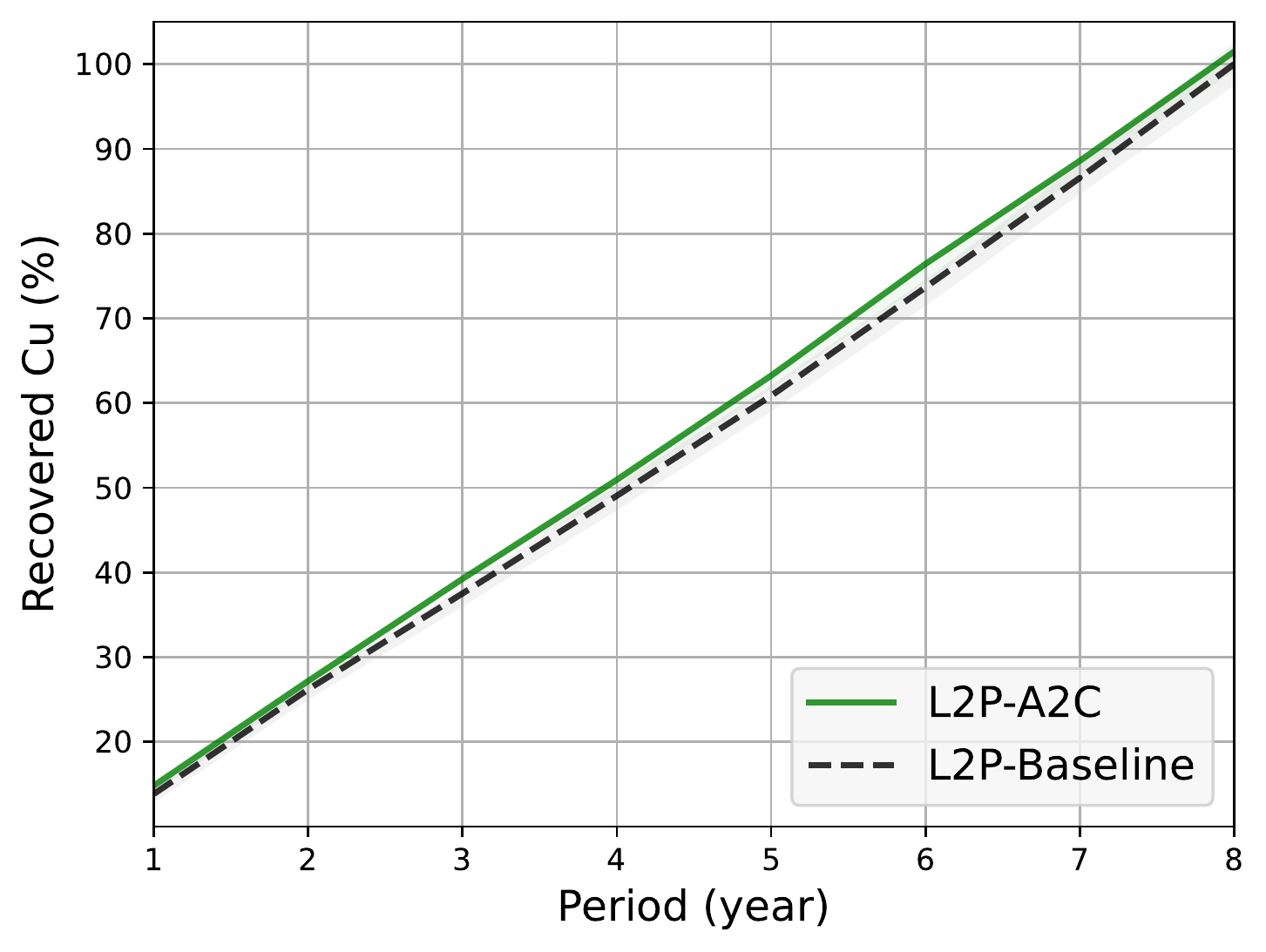}
\end{minipage}
\caption{Evolution of NPV (left plot) and recovered copper (right plot) for instance $I_4$.}
\label{fig:I4_2}
\end{figure}

To compare the solutions obtained for the instance $I_4$, Figure~\ref{fig:I4_2} shows the NPV and recovered copper for the solution obtained using L2P--A2C (the best performing variant of the L2P hyper-heuristic) and L2P--Baseline. The stopping time criterion used is equal to the execution time it took L2P--A2C to yield a solution value that is 1\% away from $Z*$. The cash flow and recovered metal for the solution obtained using L2P--A2C are 2-3\% higher than those obtained using L2P--Baseline.

\subsubsection{Third stage of testing}

\begin{figure}
\begin{minipage}[t]{.4\textwidth}
\centering
\includegraphics[width=\textwidth]{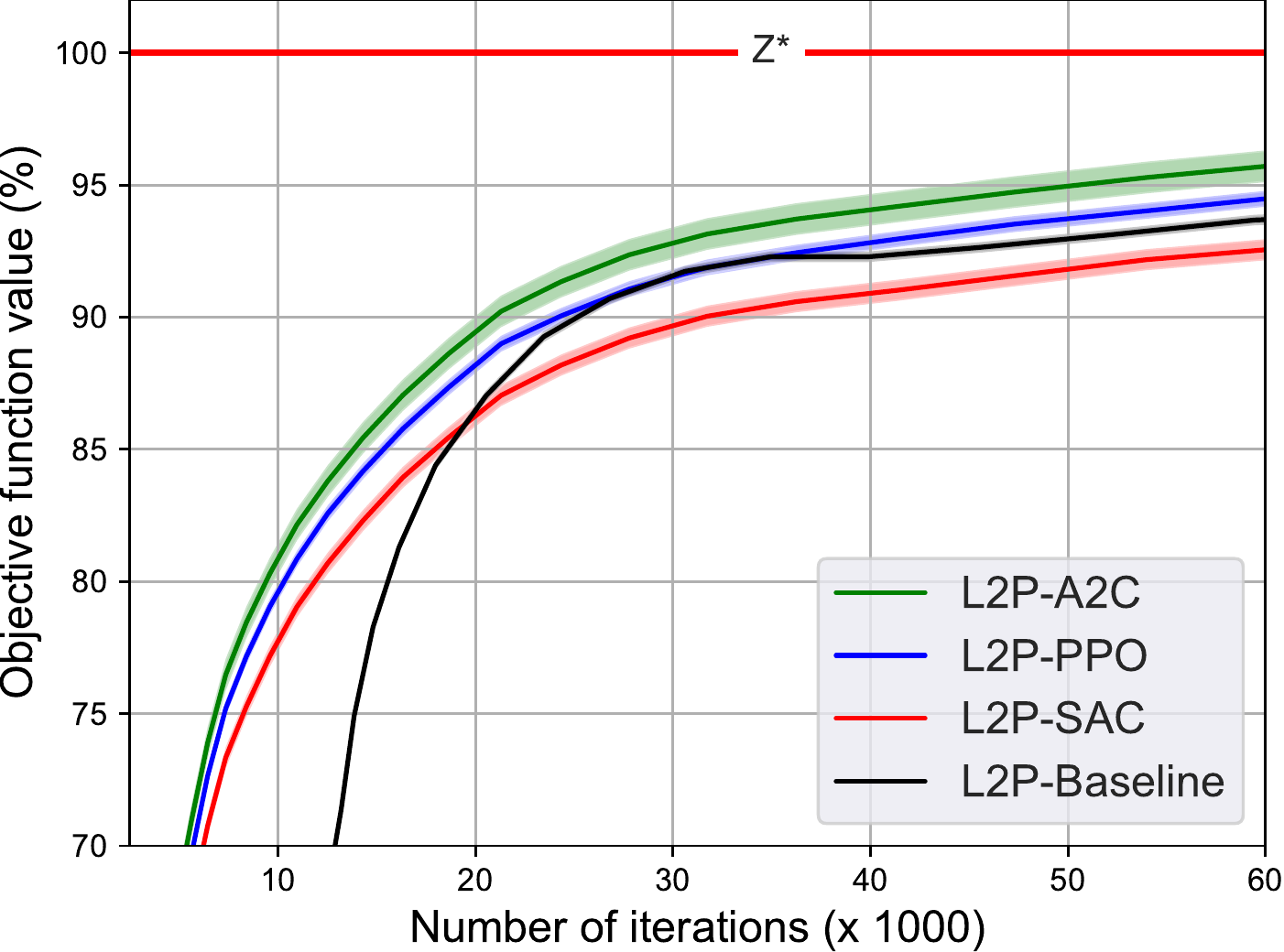}
\end{minipage}
\hfill
\begin{minipage}[t]{.4\textwidth}
\centering
\includegraphics[width=\textwidth]{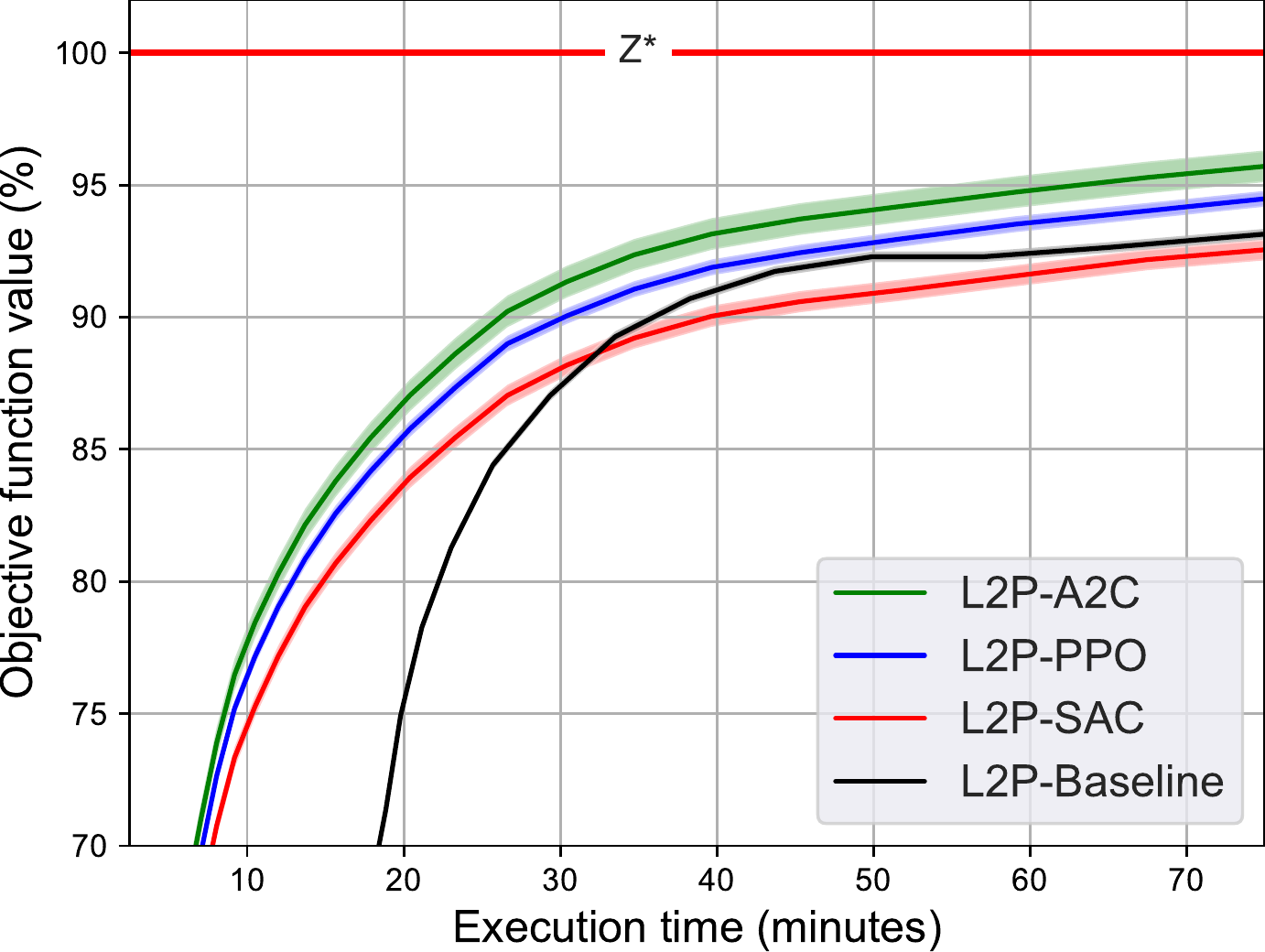}
\end{minipage}  
\caption{Evolution of the objective function value (in \%) for Instance $I_5$ with respect to the number of iterations ($\times1000$) (left plot) and the execution time (in minutes) (right plot).
}
\label{fig:day3}
\end{figure}

As in Figure~\ref{fig:day3}, L2P--A2C outperforms all other variants, reducing the number of iterations and the execution time by $45$--$50\%$, compared to L2P--Baseline. L2P--PPO reduces the number of iterations and the execution time by only $15$--$20\%$.

Note that during optimization, L2P--Baseline outperforms L2P--SAC by providing a better solution. This is only observed when comparing L2P--Baseline and L2P--SAC and can be explained by the fact that the SAC policy is non-optimal and minimally differentiates states during optimization. Although such non-optimality could be addressed by pre-training the L2P--SAC on $I_5$ instead of using pre-trained weights, it is crucial for users that the pre-trained solver be as generalizable as possible so that no tuning is required. Thus, using L2P--SAC in a real-life scenario is impractical. In contrast, L2P--A2C and L2P--PPO emerge as better variants of the L2P hyper-heuristic.

\begin{figure}
\begin{minipage}[t]{.45\textwidth}
\centering
\includegraphics[width=\textwidth]{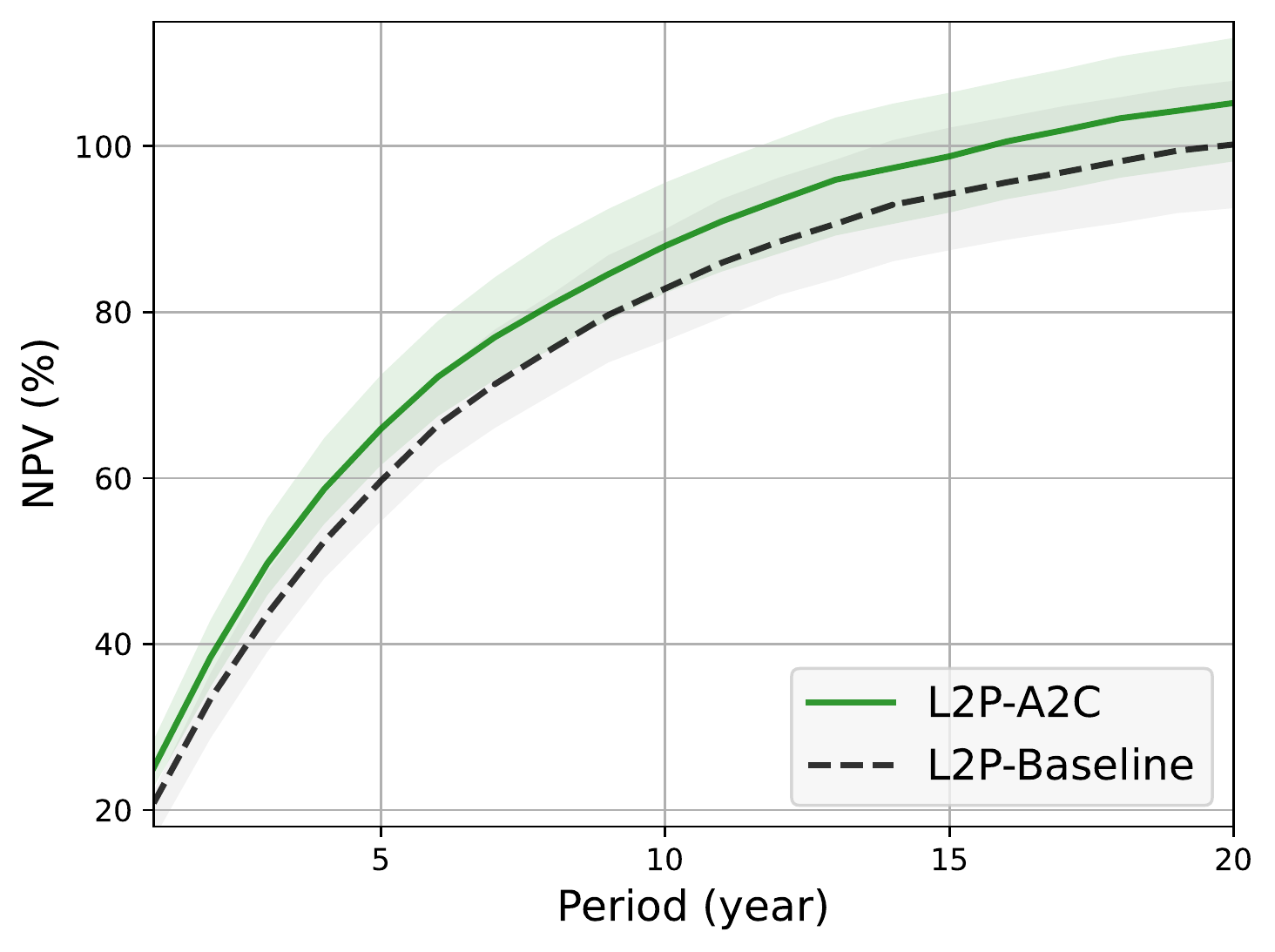}
\end{minipage}
\hfill
\begin{minipage}[t]{.45\textwidth}
\centering
\includegraphics[width=\textwidth]{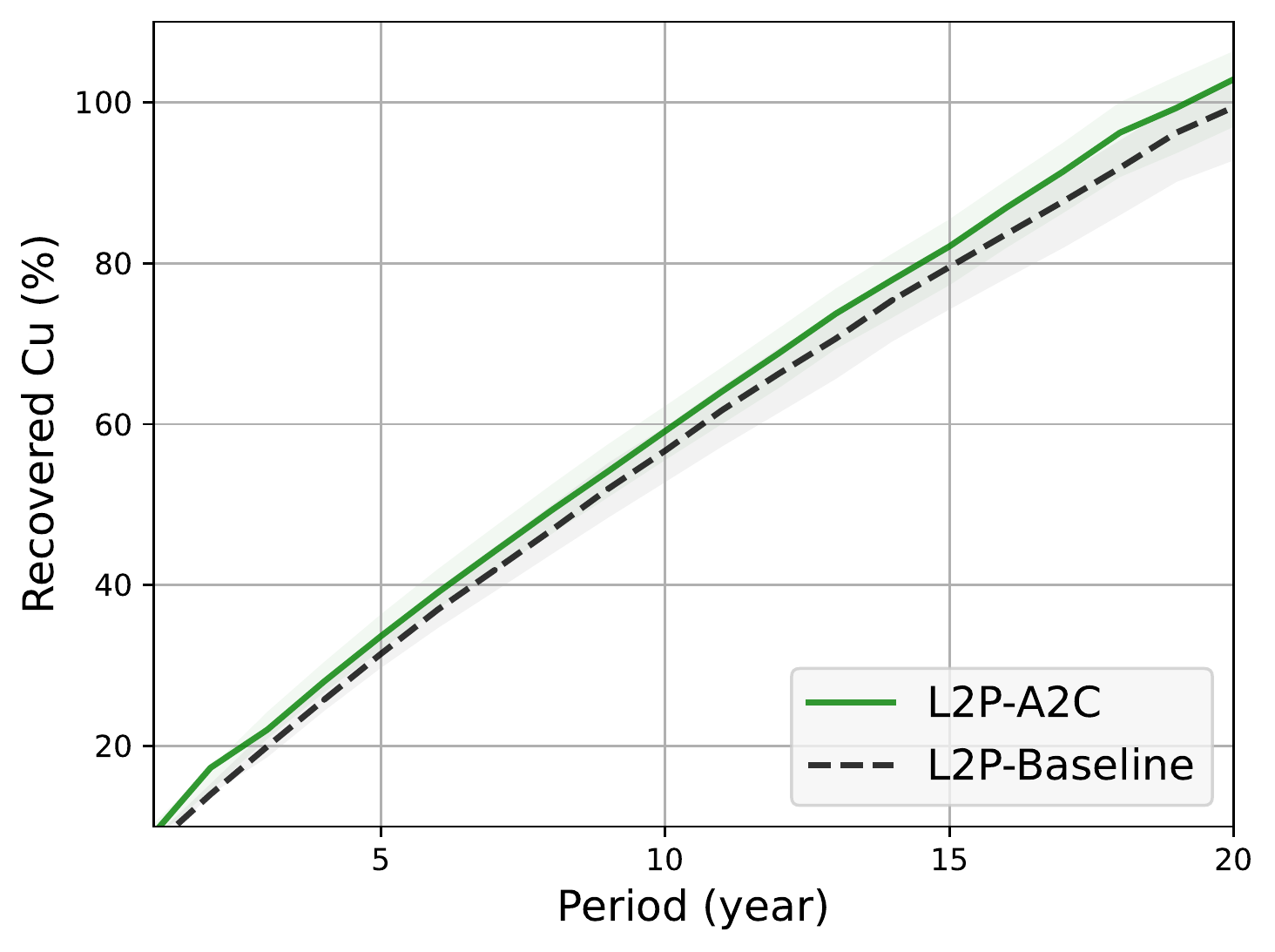}
\end{minipage}
\caption{Evolution of NPV (left plot) and recovered copper (right plot) for instance $I_5$.}
\label{fig:I5_2}
\end{figure}

To compare the solutions obtained for the instance $I_5$, Figure~\ref{fig:I5_2} shows the NPV and recovered copper for the solution obtained using L2P--A2C (the best performing variant of the L2P hyper-heuristic) and L2P--Baseline. The stopping time criterion used is equal to the execution time it took L2P--A2C to yield a solution value that is 1\% away from $Z*$. The cash flow and recovered metal (for the solution obtained using L2P--A2C) are $2-3\%$ higher than those obtained using L2P--Baseline.

\section{Conclusions}\label{sec:conclusion}

This paper presents a self-managed hyper-heuristic based on multi-neighborhood simulated annealing with adaptive neighborhood search, where state-of-the-art reinforcement learning agents are used to adapt the search and guide it towards better solutions. By defining a neighborhood structure, the reinforcement learning component learns to guide the search through the solution space landscape during the optimization process, while sampling towards the most promising heuristics.
We apply the proposed L2P hyper-heuristic to the simultaneous stochastic optimization of mining complexes (SSOMC), a large-scale stochastic combinatorial optimization problem of particular interest to the mining industry. Results comparing different variants of the L2P hyper-heuristic (with and without reinforcement learning) show its effectiveness on several case studies, reducing the number of iterations by $30$--$50\%$ and the computational time by $30$--$45\%$, compared to the baseline using a standard score function.
Since the SSOMC is similar to other scheduling problems, exploring this research direction to find more efficient methods to solve such problems becomes even more justified. Indeed, unlike problem-specific heuristics, the presented solution method can be generalized to other problems as in airlines (airline crew scheduling), bus (shift scheduling), trucking (routing problems), and rail industries. The only requirement would be to implement a set of low-level heuristics, none of which would need to be well-tuned.

\section*{Acknowledgments}
The work in this paper was funded by the National Sciences and Engineering Research Council (NSERC) of Canada CRD Grant 500414-16, the NSERC Discovery Grant 239019, and the industry consortium members of COSMO Stochastic Mine Planning Laboratory (AngloGold Ashanti, BHP, De Beers/AngloAmerican, IAMGOLD Corporation, Kinross Gold, Newmont Mining Corporation, and Vale), and the Canada Research Chairs Program.

\bibliography{Bibliography}

\newpage

\appendix

\section{Problem formulation of the two-stage stochastic integer program}\label{app:formulation}

In what follows, the mathematical formulation of the simultaneous stochastic optimization of mining complexes (SSOMC)~\citep{goodfellow2016,goodfellow2017} is outlined. All sets and indices used in the formulation are presented in Table~\ref{tab:indices}. All inputs and parameters are presented in Table~\ref{tab:inputs_parameters}.

\begin{table}[H]
\renewcommand{\arraystretch}{1.15}
\centering
\small
\caption{List of sets and indices}\label{tab:indices}
\vspace{2mm}
\begin{tabular}{|l|l|}
\hline
\textbf{Sets/Indices} & \textbf{Description}\\
\hline
$p \in \mathbb{P}$ & primary attributes, sent  from one location in the mining complex to another \\
& (e.g.,  metal tonnages, or total tonnage) \\
\hline
$h \in \mathbb{H}$ & hereditary attributes, variables of interest at specific locations (nodes) in the\\
& mining complex, but not necessarily forwarded between locations\\
& (e.g., mining, stockpiling and processing costs) \\
\hline
$i \in \mathcal{S} \cup \mathcal{D} \cup \mathcal{M}$ & locations, within the set of mines $\mathcal{M}$, stockpiles $\mathcal{S}$, and other destinations $\mathcal{D}$\\
\hline
$b \in \mathbb{B}_{m}$ & block $b$ within the set of blocks $\mathbb{B}_{m}$ at mine $m \in \mathcal{M}$ \\
\hline
$s \in \mathbb{S}$ & a set of equally probable scenarios (simulations) describing the combined \\
& simulations from different sources of uncertainty (e.g., geological uncertainty) \\
\hline
$t \in \mathbb{T}$ & time period within the life-of-mine\\
\hline
$g \in \mathbb{G}$ & a group $\mathbb{G}$ is generated as a pre-processing step using clustering to reduce\\
& the number of policy decisions variables based on multiple elements within\\
& each material type.\\ 
\hline
\end{tabular}
\end{table}

\begin{table}[H]
\renewcommand{\arraystretch}{1.15}
\centering
\small
\caption{List of inputs and parameters}\label{tab:inputs_parameters}
\vspace{2mm}
\begin{tabular}{|l|l|}
\hline
\textbf{Input/Parameter} & \textbf{Description}\\
\hline
$\beta_{p, b, s}$ & Simulated block attributes. \\
\hline
$\mathbb{O}(b)$ &  Block extraction precedence relationships. \\
\hline
$\theta_{b, g, s}$ & Block group memberships \\
\hline
$\mathcal{O}(i)$ and $\mathcal{J}(i),~\forall i \in \mathcal{S} \cup \mathcal{D} \cup \mathbb{G}$ & A model of the mining complex, where $\mathcal{O}(i)$ is the set of locations\\
& that receive materials from location $i$ and $\mathcal{J}(i)$ is the set of\\
& locations in the mining complex that send material to a location $i$.
\\
\hline
$f_{h}(p, i)$ & A model of the hereditary attribute transformation functions\\
\hline
$p_{h, i, t}$ & Time-discounted price (or cost) per unit of attribute (typically a \\
& discount rate that is used to calculate the NPV) \\
\hline
$U_{h, i, t}$ and $L_{h, i, t}$ & Upper- and lower-bounds for each attribute (typically tonnage,\\
& metal production and product quality constraints).\\
\hline
$c_{h, i, t}^{+}$and $c_{h, i, t}^{-}$ & Penalty costs used to penalize deviations from the upper and\\
& lower-bounds. These penalty costs may be time-varied to provide \\
& a geological risk discounting, i.e. $c_{h, i, t}=c_{h, i} /\left(1+grd_{h, i}\right)^{t}$, where $c_{h, i}$\\
& is a base penalty cost and $grd_{h, i}$ is the geological risk discount rate\\
& for the attribute of interest $h$. \\
\hline
\end{tabular}
\end{table}

The objective function mathematically formulated in Eq.~\eqref{Eq:SSOMC} maximizes the expected NPV of the mining complex and minimizes the expected recourse costs incurred whenever the stochastic constraints are violated.
To optimize said objective function, three types of decision variables are defined: (i) production scheduling decisions ($x_{b, t}$) defining whether a block is extracted per period, (ii) destination policy decisions ($z_{g, j, t}$) defining whether a group of material is sent to a destination per period, (iii) processing stream decisions ($y_{i, j, t, s}$) defining the proportion of product sent from a location to a destination.
Since the SSOMC is a constrainted optimization problem, it involves a set of constraints that are formulated in Eq.~\eqref{Eq:Form2}-\eqref{Eq:Form9}. Constraints~\eqref{Eq:Form2}-\eqref{Eq:Form3} are scenario-independent where Eq.~\eqref{Eq:Form2} guarantee that each block $i$ is mined at most once during the horizon and enforce mining precedence, and Eq.~\eqref{Eq:Form3} ensure that a group of material is sent to a single destination. Constraints~\eqref{Eq:Form4}-\eqref{Eq:Form8} are scenario-dependent and ensure the proper calculation of attributes throughout the mining complex. In particular, Eq.~\eqref{Eq:Form4} calculates the quantity of primary attributes at each location and ensue mass balancing (processing stream calculations). Eq.~\eqref{Eq:Form5} calculates the values of the hereditary attributes based on the values of the primary attributes (attribute calculations). Eq.~\eqref{Eq:Form6}-\eqref{Eq:Form7} calculate the amount of constraint violation from upper and lower bounds imposed on hereditary attributes (deviation constraints). Finally, Eq.~\eqref{Eq:Form8} is optional and calculates the end-of-year stockpile quantity.

\textit{Objective Function:}
\begin{equation}\label{Eq:SSOMC}
\max \frac{1}{|\mathbb{S}|} \underbrace{\sum_{s \in \mathbb{S}} \sum_{t \in \mathbb{T}} \sum_{h \in \mathbb{H}} p_{h, i, t} \cdot v_{h, i, t, s}}_{\text {Discounted costs and revenues }}-\frac{1}{|\mathbb{S}|} \underbrace{\sum_{s \in \mathbb{S}} \sum_{t \in \mathbb{T}} \sum_{h \in \mathbb{H}} c_{h, i, t}^{+} \cdot d_{h, i, t, s}^{+}+c_{h, i, t}^{-} \cdot d_{h, i, t, s}^{-}}_{\text {Penalties for deviations from targets }}.    
\end{equation}

\textit{Subject to:}

\begin{equation}
\begin{aligned}\label{Eq:Form2}
\sum_{t \in \mathbb{T}} x_{b, t} & \leq 1 \quad \forall b \in \mathbb{B}_{m}, m \in \mathcal{M}, \\
x_{b, t} & \leq \sum_{t^{\prime}=1}^{t} x_{u, t^{\prime}} \quad \forall b \in \mathbb{B}_{m}, m \in \mathcal{M}, u \in \mathbb{O}(b), t \in \mathbb{T} .
\end{aligned}
\end{equation}

\begin{equation}\label{Eq:Form3}
\sum_{j \in \mathcal{O}(g)} z_{g, j, t}=1 \quad \forall g \in \mathbb{G}, t \in \mathbb{T}.   
\end{equation}

\begin{align}\label{Eq:Form4}
v_{p, j,(t+1), s}&=\underbrace{v_{p, j, t, s} \cdot\left(1-\sum_{k \in \mathcal{O}(j)} y_{j, k, t, s}\right)}_{\text {Material from previous period }}+\underbrace{\sum_{i \in \mathcal{J}(j) \backslash \mathbb{G}} r_{p, i, t, s} \cdot v_{p, i, t, s} \cdot y_{i, j, t, s}}_{\text {Incoming materials from other locations }} \nonumber \\
&\qquad+\underbrace{\sum_{g \in \mathcal{J}(j) \cap \mathbb{G}}\left(\sum_{b \in \mathbb{B}_{m}} \sum_{m \in \mathcal{M}} \theta_{b, g, s} \cdot \beta_{p, b, s} \cdot x_{b,(t+1)}\right) \cdot z_{g, j,(t+1)}}_{\text {Materials sent directly from mines }} \nonumber \\
&\qquad \forall p \in \mathbb{P}, j \in \mathcal{S} \cup \mathcal{D}, t \in \mathbb{T}, s \in \mathbb{S} \nonumber\\
&\qquad \sum_{j \in \mathcal{O}(i)} y_{i, j, t, s}=1 \quad \forall i \in \mathcal{D}, t \in \mathbb{T}, s \in \mathbb{S} \nonumber\\
&\qquad \sum_{j \in \mathcal{O}(i)} y_{i, j, t, s} \leq 1 \quad \forall i \in \mathcal{S}, t \in \mathbb{T}, s \in \mathbb{S}.
\end{align}

\begin{align}\label{Eq:Form5}
v_{h, i, t, s} &=f_{h}(p, i) \quad \forall h \in \mathbb{H}, i \in \mathcal{S} \cup \mathcal{D} \cup \mathcal{M}, t \in \mathbb{T}, s \in \mathbb{S} \nonumber \\
v_{p, m, t, s} &=\sum_{b \in \mathbb{B}_{m}} \beta_{p, b, s} \cdot x_{b, t} \quad \forall m \in \mathcal{M}, p \in \mathbb{P}, t \in \mathbb{T}, s \in \mathbb{S} .
\end{align}

\begin{align}\label{Eq:Form6}
&v_{h, i, t, s}-d_{h, i, t, s}^{+} \leq U_{h, i, t} \quad \forall h \in \mathbb{H}, t \in \mathbb{T}, s \in \mathbb{S} \nonumber \\
&v_{h, i, t, s}+d_{h, i, t, s}^{-} \geq L_{h, i, t} \quad \forall h \in \mathbb{H}, t \in \mathbb{T}, s \in \mathbb{S} .
\end{align}


\begin{align}\label{Eq:Form7}
&r_{p, i, t, s}=1 \quad \forall p \in \mathbb{P}, i \in \mathcal{S}, t \in \mathbb{T}, s \in \mathbb{S} \nonumber \\
&r_{p, i, t, s}=v_{h, i, t, s} \quad \forall p \in \mathbb{P}, i \in \mathcal{D}, t \in \mathbb{T}, s \in \mathbb{S} .
\end{align}

\begin{equation}\label{Eq:Form8}
v_{h, i, t, s}=v_{p, i, t, s} \cdot\left(1-\sum_{i \in \mathcal{O}(i)} y_{i, j, t, s}\right) \quad \forall h \in \mathbb{H}, i \in \mathcal{S}, t \in \mathbb{T}, s \in \mathbb{S}.
\end{equation}

\begin{align}\label{Eq:Form9}
&x_{b, t} \in\{0,1\} \quad \forall b \in \mathbb{B}_{m}, m \in \mathcal{M}, t \in \mathbb{T} \nonumber \\
&z_{g, j, t} \in\{0,1\} \quad \forall g \in \mathbb{G}, j \in \mathcal{O}(g), t \in \mathbb{T} \nonumber \\
&y_{i, j, t, s} \in[0,1] \quad \forall i \in \mathcal{S} \cup \mathcal{D}, j \in \mathcal{O}(i), t \in \mathbb{T}, s \in \mathbb{S} \nonumber\\
&v_{p, i, t, s} \geq 0 \quad \forall p \in \mathbb{P}, i \in \mathcal{S} \cup \mathcal{D} \cup \mathcal{M}, t \in \mathbb{T}, s \in \mathbb{S} \nonumber\\
&v_{h, i, t, s} \in \mathbb{R} \quad \forall h \in \mathbb{H}, i \in \mathcal{S} \cup \mathcal{D} \cup \mathcal{M}, t \in \mathbb{T}, s \in \mathbb{S} \nonumber\\
&r_{p, i, t, s} \in[0,1] \quad \forall p \in \mathbb{P}, i \in \mathcal{S} \cup \mathcal{D}, t \in \mathbb{T}, s \in \mathbb{S} \nonumber\\
&d_{h, i, t, s}^{+}, d_{h, i, t, s}^{-} \geq 0 \quad \forall h \in \mathbb{H}, i \in \mathcal{S} \cup \mathcal{D} \cup \mathcal{M}, t \in \mathbb{T}, s \in \mathbb{S} .
\end{align}

\newpage
\section{Reinforcement learning methods}\label{app:rl}

In what follows, the pseudo-codes for the employed reinforcement learning agents are detailed, namely the Advantage Actor-Critic (A2C)~\citep{a2c} in Algorithm~\ref{algo:1}, the Soft Actor-Critic (SAC)~\citep{sac} in Algorithm~\ref{algo:2}, and the Proximal Policy Optimization (PPO)~\citep{ppo} in Algorithm~\ref{algo:3}. The notation that are used used are reported in Table~\ref{policy_gradient_Notations}.

\begin{algorithm}[H]
	\caption{Advantage Actor-Critic (A2C)}
	\label{algo:1}
	\begin{algorithmic}[1]
		\Procedure {A2C}{$\theta,w$}
		\State let $\theta, w$ be global parameters, and $\theta', w'$ are thread-specific parameters;
		\State $t\gets 1$;
		\While{$T \leq T_{\MAX}$}
		\State Reset gradient $\phi_\theta \gets0, \; \phi_{w} \gets0$, and set $\theta'\gets\theta,\; w'\gets w$;
		\State let $t_{\start} \gets t$ and a starting state $s_t$ is sampled;
		\While {$s_t \neq \TERMINAL \; \& \; t - t_{\start}\leq t_{\MAX}$}
		\State Pick the action $a \sim \pi_{\theta}(a_t,s_t)$ to receive the reward $r_t$ and move to next state $s_{t+1}$;
		\State $T \gets T+1, \; t \gets t+1$;
		\EndWhile
		\If{$s_t = \TERMINAL$ }
		\State $R \gets 0$;
		\Else
		\State $R \gets V_{w}(s_t)$;
		\EndIf
		
		\For{$i \gets (t-1)$ \textbf{to}  $t_{\start}$}
		\State $R \gets \gamma R + R_i$;
		\State Accumulate gradients w.r.t $\theta': \phi\theta \gets \phi\theta + \nabla_{\theta'} \log \pi_{\theta'}(a_i \vert s_i)(R - V_{w'}(s_i))$;
		\State Accumulate gradients w.r.t $w': \phi w \gets \phi w + 2 (R - V_{w'}(s_i)) \nabla_{w'} (R - V_{w'}(s_i))$
		\EndFor
		\State $\theta \gets \phi_\theta, \; w \gets \phi_{w}$
		\EndWhile
		\EndProcedure
		
	\end{algorithmic}
\end{algorithm}

\begin{algorithm}[H]
	\caption{Soft Actor-Critic (SAC)}
	\label{algo:2}
	\begin{algorithmic}[1]
	\Require The learning rates $\lambda_\pi, \lambda_Q, \lambda_V$ for $\pi_{\theta}, Q_{w}, V_\psi$; The $\imath$ as weighting factor for moving average.
	\Procedure {SAC}{$\theta,w, \psi, \bar{\psi}$}
	\State initialize the parameters $\theta,w, \psi, \bar{\psi}$.
	\ForAll{iteration}
\ForAll{Environment Setup}
	\State $a_t \sim \pi_{\theta}(a_t|s_t)$, and $s_{t+1} \sim \rho{\pi}(s_{t+1}|a_t,s_t)$
	\State $\mathcal{Z}\gets \mathcal{Z}\cup{(s_t,a_t,r((s_t,a_t)),s_{t+1}}$
	\EndFor
	\ForAll{Gradient Update Step}
	\State $\psi\gets\psi-\lambda_V\nabla_\psi Y_V(\psi)$, and $w\gets w-\lambda_Q\nabla_{w} Y_Q(w)$
	\State $\theta\gets\theta-\lambda_\pi\nabla_\theta Y_\pi(\theta)$, and $\bar{\psi}\gets\imath\psi+((1-\imath)\bar{\psi})$
	\EndFor
	\EndFor
	\EndProcedure
\end{algorithmic}
\end{algorithm}

 \begin{algorithm} [H]
 	\caption{Proximal Policy Optimization (PPO)}
 	\label{algo:3}
 	\begin{algorithmic}[1]
 		\Require The initial policy parameters $\theta$, the clipping threshold $\epsilon$.
 		\Procedure {PPO}{$\theta,\epsilon$}
 		\For{$i\gets 1$, \ldots}.
 		\State Collect the set of partial trajectories $\mathcal{D}$ for policy $\pi_i=\pi(\theta_i)$
 		\State Estimate Advantages $\hat{\mathcal{W}}^{\theta_i}_t$ by utilizing some estimation algorithm.
 		\State Calculate the policy update $\theta_{i+1}=\argmax Y^{\CLIP}_{\theta_i}$
 		\State 	$Y^{\CLIP}_{\theta_i} (\theta) = \mathbb{E}_{t\sim{\theta_i}} \left[ \sum_{t=0}^{t^{\max}} \left[ \min( r(\theta) \hat{\mathcal{W}}_{\theta_\text{old}}(s, a), {\CLIP}(r(\theta), 1 - \epsilon, 1 + \epsilon) \hat{\mathcal{W}}_{\theta_{\old}}(s, a))\right]\right]$ 
 		
 		\EndFor
 		\EndProcedure
 		
 	\end{algorithmic}
 \end{algorithm}

\end{document}